\documentclass[leqno,12pt]{article}  %putting formula numbers on the left side
\usepackage{amsmath}
\usepackage{amssymb}

\usepackage[german,english]{babel}
\usepackage{amsthm}
\usepackage{xypic}  % commutative diagrams

\author{\textsc{Elmar Grosse-Kl\"onne}}
\date{}

\theoremstyle{plain} 
\newtheorem{satz}{Theorem}[section]  %@@
\newtheorem{kor}[satz]{Corollary}  %@@
\newtheorem{lem}[satz]{Lemma}  %@@
\newtheorem{pro}[satz]{Proposition}  %@@
\theoremstyle{remark}

\theoremstyle{definition}

\newcommand{\0}{\ensuremath{\overrightarrow{0}}}

\begin{document}

% \addtolength{\textwidth}{1.0in}
% \setlength{\hoffset}{-.5in}
%
%\addtolength{\topmargin}{-34pt}
%\addtolength{\textheight}{68pt}
%\baselineskip13,6pt
%\setlength{\parindent}{0pt}
%\setlength{\parskip}{0.8ex plus 0.2ex minus 0.2ex}
%\setlength{\oddsidemargin}{3cm}
%\frenchspacing
%\sloppy
%\pagestyle{myheadings}
%\markboth{}{}

%\textwidth14cm
%\textheight9in
%\hoffset+1in
%\voffset-0,8cm
\title{ Locally algebraic automorphisms of the ${\rm PGL}_2(F)$-tree and ${\mathfrak o}$-torsion representations }

\begin{center}{\bf Locally algebraic automorphisms of the ${\rm
      PGL}_2(F)$-tree
    and ${\mathfrak o}$-torsion representations }\\Elmar Grosse-Kl\"onne
\end{center}

\begin{abstract} For a local field $F$ and an Artinian local coefficient ring $\Lambda$ with the same positive residue characteristic $p$ we define, for any $e\in{\mathbb N}$, a category ${\mathfrak C}^{(e)}(\Lambda)$ of ${\rm GL}_2(F)$-equivariant coefficient systems on the Bruhat-Tits tree $X$ of ${\rm PGL}_2(F)$. There is an obvious functor from the category of ${\rm GL}_2(F)$-representations over $\Lambda$ to ${\mathfrak C}^{(e)}(\Lambda)$. If $F={\mathbb Q}_p$ then ${\mathfrak C}^{(1)}(\Lambda)$ is equivalent to the category of smooth ${\rm GL}_2({\mathbb Q}_p)$-representations over $\Lambda$ generated by their invariants under a pro-$p$-Iwahori subgroup. For general $F$ and $e$ we show that the subcategory of all objects in ${\mathfrak C}^{(e)}(\Lambda)$ with trivial central character is equivalent to a category of representations of a certain subgroup of ${\rm Aut}(X)$ consisting of "locally algebraic
automorphisms of level $e$". For $e=1$ there is a functor from this
category to that of modules over the (usual) pro-$p$-Iwahori Hecke
algebra; it is a bijection between irreducible objects.

Finally, we present a parallel of Colmez' functor $V\mapsto {\bf
  D}(V)$: to objects in ${\mathfrak C}^{(e)}(\Lambda)$ (for {\it
  any} $F$) we assign certain \'{e}tale $(\varphi,\Gamma)$-modules over an Iwasawa algebra ${\mathfrak
  o}[[\widehat{N}^{(1)}_{0,1}]]$ which contains the (usually considered) Iwasawa algebra ${\mathfrak
  o}[[{N}_{0}]]$. This assignment preserves finite generation.\\


{\it AMS 2010} Mathematics subject classification: 22E50, 20G25, 20E42 
\end{abstract}

\section{Introduction}

Let $F$ be a local field with residue characteristic $p>0$ and uniformizer
$p_F\in{\mathcal O}_F$ generating the maximal ideal ${\mathfrak p}_F$ of the ring of integers ${\mathcal O}_F$. Let $G$ be the group of $F$-rational
points of a reductive algebraic group over $F$. An important tool in the smooth representation
theory of $G$ on vector spaces over the complex numbers ${\mathbb C}$ is the localization
technique which has been systematically developed by Schneider and Stuhler in
their work
\cite{ss}. Assigning to a smooth (admissible, finite length) $G$-representation on a ${\mathbb C}$-vector
space $V$ and a simplex
$\tau$ in the
Bruhat-Tits building of $G$ the space of invariants of $V$ under a suitable
open subgroup of $G$ fixing $\tau$, one obtains a $G$-equivariant (homological) coefficient
system ${\mathcal F}_V$ on $X$. If this
assignment is carried out with appropriate care then $V$ can be recovered from
${\mathcal F}_V$ as $V=H_0(X,{\mathcal F}_V)$, and in this way, the study of smooth (admissible, finite length) complex
$G$-representations is transformed into the study of coefficient
systems on $X$ --- in a sense these coefficient
systems are 'smaller' objects, accessible by the representation theory of
finite groups. These constructions work well also for smooth $G$-representations on vector
spaces over fields of positive characteristic different from $p$.

On the other hand, if we ask for smooth $G$-representations on vector
spaces $V$ over a field $k$ of characteristic $p$, then analogous assignments
$V\mapsto{\mathcal F}_V$ are much weaker in general; typically, they do not
allow to recover $V$. There seems to be basically only one example class of
smooth (and possibly supercuspidal/supersingular) $G$-representations over $k$ for which the classical
(complex) theory carries over to wide extent: this is the case where $G={\rm
  GL}_2({\mathbb Q}_p)$ (or $G={\rm
  SL}_2({\mathbb Q}_p)$, or $G={\rm
  PGL}_2({\mathbb Q}_p)$) and
where the smooth $G$-representations considered are generated by their invariants
under a pro-$p$-Iwahori subgroup $U^{(1)}_{\sigma}$ of $G$. Namely, the
category of these smooth $G$-representations is equivalent to a category of
$G$-equivariant coefficient systems on the Bruhat-Tits tree $X$ of ${\rm
  PGL}_2({\mathbb Q}_p)$ satisfying a simple and natural
axiomatic (i.e. the category  ${\mathfrak C}^{(1)}(k)$ below). See \cite{jalg} or Theorem \ref{jalgcentral} below for the precise
statement.

The purpose of this paper is to discuss similar concepts for the groups $G={\rm
  GL}_2(F)$ for general $F$. As in \cite{ss} we fix an index $e\ge1$
(the 'level') and, for an edge $\eta$ of the Bruhat-Tits tree $X$ of ${\rm
  PGL}_2(F)$, we consider the open subgroup $U^{(e)}_{\eta}$ of $G$ 'of level $e$'
which fixes $\eta$. We fix an edge $\sigma$ of $X$. For a ring $\Lambda$ let ${\mathfrak C}^{(e)}(\Lambda)$ denote the category of
$G$-equivariant homological coefficient systems of $\Lambda$-modules ${\mathcal F}$ on $X$ such
that for any vertex $x$ and any edge $\eta$ with $x\in\eta$ the transition
map ${\mathcal
  F}(\eta)\to {\mathcal
  F}(x)$ is injective, its image is ${\mathcal
  F}(x)^{U^{(e)}_{\eta}}$ and generates ${\mathcal
  F}(x)$ as a representation of the stabilizer of $x$ in $G$. There is an obvious functor $V\mapsto {\mathcal F}_V^{(e)}$ from the category of $G$-representation on $\Lambda$-modules to the category ${\mathfrak C}^{(e)}$; it satisfies ${\mathcal F}_V^{(e)}(\sigma)=V^{U^{(e)}_{\sigma}}$. If $\Lambda={\mathbb C}$ then the category of smooth, admissible,
finite length $G$-representations over ${\mathbb C}$ generated by
$V^{U^{(e)}_{\sigma}}$ embeds into
a full subcategory of
${\mathfrak C}^{(e)}({\mathbb C})$ by means of this functor $V\mapsto {\mathcal F}_V^{(e)}$. Therefore it is natural to ask for the relevance of the category ${\mathfrak C}^{(e)}(\Lambda)$ for arbitrary $\Lambda$. Is it equivalent to a suitable category of $G$-representations ?

Let ${\mathfrak C}_0^{(e)}(\Lambda)$ denote the
subcategory of all ${\mathcal F}\in {\mathfrak C}^{(e)}(\Lambda)$ on which the
action of $G$ factors through ${\rm PGL}_2(F)$. The basic
observation of the present paper is that the ${\rm PGL}_2(F)$-action on any ${\mathcal
  F}\in{\mathfrak C}_0^{(e)}(\Lambda)$ and also on its homology $H_0(X,{\mathcal
  F})$ naturally extends to a much larger group $\widehat{G}^{(e)}$ containing
${\rm PGL}_2(F)$ and contained in the automorphism group ${\rm Aut}(X)$ of the
tree $X$. Briefly, an element of ${\rm Aut}(X)$ belongs to $\widehat{G}^{(e)}$ if and
only if for any edge $\eta$ of $X$ it acts on the ball of radius
$e+\frac{1}{2}$ around $\eta$ like an element of ${\rm
  PGL}_2(F)$. We call the elements of $\widehat{G}^{(e)}$ locally algebraic
automorphisms of $X$ of level $e$.\footnote{In particular we see that the smooth, admissible,
finite length ${\rm PGL}_2(F)$-representations $V$ over ${\mathbb C}$
generated by $V^{U^{(e)}_{\sigma}}$ automatically carry an action by the
larger group $\widehat{G}^{(e)}$. For $F\ne{\mathbb Q}_p$ this fails if ${\mathbb C}$ is replaced by
$k$ (as follows e.g. from \cite{brpa}), and from
the point of view of the present paper, the failure of this
principle is the reason for, or the manifestation of, the difference between
the smooth
${\rm PGL}_2(F)$-representation theory over ${\mathbb C}$ and over $k$.} 

The subgroups $U^{(e)}_{\eta}$ of $G$ have as natural analogs certain pro-$p$-subgroups 
$\widehat{U}^{(e)}_{\eta}$ of $\widehat{G}^{(e)}$. Given a
$\widehat{G}^{(e)}$-representation $V$ we assign to it the coefficient system
$\widehat{\mathcal F}^{(e)}_V\in {\mathfrak C}_0^{(e)}(\Lambda)$ with $\widehat{\mathcal
  F}_V^{(e)}(\eta)=V^{\widehat{U}^{(e)}_{\eta}}$ for edges $\eta$. Let
$\widehat{\mathfrak R}^{(e)}_0(\Lambda)$ denote the category of
$\widehat{G}^{(e)}$-representations $V$ over $\Lambda$ generated by
$V^{\widehat{U}^{(e)}_{\sigma}}$, and smooth when regarded as representations of $\widehat{U}^{(e)}_{\sigma}$. (We find it convenient {\it not} to work with a topology on $\widehat{G}^{(e)}$; as a consequence, we do {\it not} have available the concept of a smooth $\widehat{G}^{(e)}$-representation. Instead, the subgroups $\widehat{U}^{(e)}_{\sigma}$ (and their open subgroups) mimick the role which open subgroups play in usual smooth representation theory.) Let now ${\mathfrak o}$ be a complete
discrete valuation ring with residue field $k$, and assume that $\Lambda$ is an
Artinian local ${\mathfrak o}$-algebra with residue field $k$. Then our first
main theorem is the following (Theorem \ref{centralequiv}):\\

{\bf Theorem:} {\it The assignments ${\mathcal F}\mapsto
  H_0(X,{\mathcal F})$ and $V\mapsto\widehat{\mathcal F}_V^{(e)}$ are an
  equivalence of categories between ${\mathfrak C}_0^{(e)}(\Lambda)$ and
  $\widehat{\mathfrak R}_0^{(e)}(\Lambda)$. For ${\mathcal F}\in {\mathfrak C}_0^{(e)}(\Lambda)$ and $V\in \widehat{\mathfrak R}_0^{(e)}(\Lambda)$ we have natural isomorphisms$${\mathcal
  F}\longrightarrow \widehat{\mathcal F}_{H_0(X,{\mathcal
    F})}^{(e)}\quad\quad\mbox{ and }\quad\quad H_0(X,\widehat{\mathcal
    F}_V^{(e)})\longrightarrow V.$$} \\  

Let ${\mathcal
  H}_{\Lambda}(G,{U}^{(1)}_{\sigma})\cong {\Lambda}[{U}^{(1)}_{\sigma}\backslash G/{U}^{(1)}_{\sigma}]$ denote the pro-$p$-Iwahori Hecke algebra over
$\Lambda$. For a smooth $G$-representation $V$ over $\Lambda$ the submodule $V^{{U}^{(1)}_{\sigma}}$ of
${U}^{(1)}_{\sigma}$-invariants of $V$ is naturally
an ${\mathcal H}_{\Lambda}(G,{U}^{(1)}_{\sigma})$-right module. In an early stadium of the investigations of the smooth representation theory
of $G$ over $k$ it was not clear (see e.g. \cite{vigcomp}) if,
parallel to similar results over ${\mathbb C}$, the functor $V\mapsto
V^{{U}^{(1)}_{\sigma}}$ was a bijection between isomorphism classes of smooth (admissible)
irreducible $G$-representations over $k$ with a central character, and isomorphism
classes of simple ${\mathcal H}_{k}(G,{U}^{(1)}_{\sigma})$-right modules (assuming
$k$ to be algebraically closed). For $F={\mathbb
  Q}_p$ this is indeed correct (see Vign\'{e}ras \cite{vigcomp}; this also uses important
work of Breuil). However, work of Breuil and Paskunas \cite{brpa} then showed that if $F\ne{\mathbb
  Q}_p$ there are many more smooth irreducible $G$-representations over $k$ than
expected, disproving such a correspondence if $F\ne{\mathbb
  Q}_p$. Again we suggest to look at $\widehat{G}^{(1)}$-representations
instead of $G$-representations. We observe that also for any $V\in \widehat{\mathfrak
  R}_0^{(1)}(\Lambda)$ the $\Lambda$-module $V^{\widehat{U}^{(1)}_{\sigma}}$ of $\widehat{U}^{(1)}_{\sigma}$-invariants of $V$ is naturally
an ${\mathcal H}_{\Lambda}(G,{U}^{(1)}_{\sigma})$-right module (but we do {\it not}
consider the Hecke algebra ${\Lambda}[\widehat{U}^{(1)}_{\sigma}\backslash
\widehat{G}^{(1)}/\widehat{U}^{(1)}_{\sigma}]$). We call $V\in\widehat{\mathfrak
  R}_0^{(1)}(\Lambda)$ admissible if $V^{\widehat{U}^{(1)}_{\sigma}}$ is a
finitely generated $\Lambda$-module. Our second main
theorem then reads (Corollary \ref{heckmain}):\\

{\bf Theorem:} {\it Assume that $k$ is algebraically closed. The functor $V\mapsto V^{\widehat{U}^{(1)}_{\sigma}}$ induces a bijection between the
  isomorphism classes of admissible irreducible
  $\widehat{G}^{(1)}$-representations and the isomorphism classes of simple
  ${\mathcal H}_{k}(G,{U}^{(1)}_{\sigma})$-right modules with trivial central
  character.}\\

In this connection let us mention Vign\'{e}ras' result \cite{vigcomp}: the number
of supersingular simple ${\mathcal H}_{k}(G,{U}^{(1)}_{\sigma})$-right modules
with a given action of the scalar matrix $p_F$ is exactly the number of
irreducible representations of the Weil group of $F$ of dimension $2$ over $k$ with a given
value of the determinant at $p_F$.\footnote{Initiated by Cartier, the structure theory of the full
automorphism group ${\rm Aut}(X)$ of $X$ (decomposition
theorems similar to those known for ${\rm PGL}_2(F)$) and its smooth complex representation theory have been thoroughly developed, see
e.g. \cite{choucrou}. It follows from results of \cite{choucrou} that the structure theory of the groups $\widehat{G}^{(e)}$
is parallel to that of ${\rm Aut}(X)$, see Theorem \ref{choucit}. Motivated by the present work
one might ask for their
representation theory over ${\mathfrak
  o}$.}   

In the final section we take a look, from the perspective of the present
paper, onto Colmez' functor from smooth ${\rm GL}_2({\mathbb
  Q}_p)$-representations on ${\mathfrak
  o}$-torsion modules (where $[{\mathfrak o}:{\mathbb Z}_p]<\infty$) to
$(\varphi,\Gamma)$-modules. First, for general $F$, let $N_0^{(e)}$ denote the pro-$p$-subgroup of $G$ consisting of
unipotent upper triangular matrices with off-diagonal entries in $p_F^{e-1}{\mathcal
  O}_F$. In $\widehat{G}^{(e)}$ we find (non-abelian) pro-$p$-subgroups
$\widehat{N}^{(e)}_{0}$ and $\widehat{N}^{(e)}_{0,1}$ which play an analogous
role, but which are much
larger than $N_0^{(e)}$ (now $\widehat{G}^{(e)}$ itself disappears from our discussion). $\widehat{N}_{0}^{(e)}$ is a
product of copies of $\widehat{N}^{(e)}_{0,1}$, and, if we set ${\mathfrak P}_k=\prod_{a\in{\mathcal
    O}_F/{\mathfrak p}^k_F}{\mathfrak p}^{k+e-1}_F$ for $k\ge0$, then
$\widehat{N}_{0,1}^{(e)}$ is abstractly the quotient of$$(\ldots({\mathfrak P}_k     \rtimes(\ldots(  {\mathfrak P}_2       \rtimes({\mathfrak P}_1       \rtimes{\mathfrak P}_0     ))\ldots))\ldots)$$by
the product of the images of all maps$$({\rm diag},-{\rm id}):\prod_{a\in{\mathcal O}_F/{\mathfrak p}^k_F}{\mathfrak p}^{k+e}_F\longrightarrow{\mathfrak P}_{k+1}\rtimes{\mathfrak P}_k.$$The $N_0^{(e)}$-action on any ${\mathcal F}\in {\mathfrak C}^{(e)}({\mathfrak
  o})$ naturally extends to a $\widehat{N}^{(e)}_0$-action on ${\mathcal F}$, hence also to a $\widehat{N}^{(e)}_0$-action on $H_0(X,{\mathcal F})$.

Take $e=1$ for the moment (in the text we consider general $e\ge1$, but see Lemma \ref{vorslev}). Given an ${\mathfrak o}$-torsion object ${\mathcal F}\in {\mathfrak C}^{(1)}({\mathfrak
  o})$ such that ${\mathcal F}(\sigma)$ is a finitely generated ${\mathfrak
  o}$-module, we faithfully copy Colmez' constructions, as reconsidered e.g. in
\cite{vigwam}: To ${\mathcal F}$ we functorially assign a finitely generated module $D({\mathcal F})$ over ${\mathcal L}(\widehat{N}^{(1)}_{0,1})={\mathfrak o}[[\widehat{N}^{(1)}_{0,1}]]$, endowed with an additional \'{e}tale structure (see the text for the precise definition: a morphism $D'({\mathcal F})\to D({\mathcal F})$ and compatible
semilinear actions by the monoid $\left( \begin{array}{cc} {\mathcal
      O}_F-\{0\}&0\\0&1\end{array}\right)$ resp. its inverse monoid, satisfying an \'{e}taleness requirement). 

Precomposing with $V\mapsto{\mathcal F}_V^{(1)}$ we obtain a functor which is a variation, available for arbitrary $F$, on Colmez' functor $V\mapsto {\bf
  D}(V)$. 

 However, if $F\ne{\mathbb Q}_p$ we are {\it not} able to assign {\it finitely generated}
  \'{e}tale $(\varphi,\Gamma)$-modules over Fontaine's
ring ${\mathcal O}_{{\mathcal E}}$ to ${\mathcal F}\in {\mathfrak C}^{(1)}({\mathfrak
  o})$.

{\it Acknowledgements:} I thank Gergely Zabradi for an important comment on an earlier draft. I thank the referee for his very careful reading of the text and for numerous detailed suggestions for improving the exposition. A substantial amount of these results has been obtained while the Deutsche Forschungsgemeinschaft (DFG) supported my position as a Heisenberg Professor at the Humboldt University at Berlin.

\section{Coefficient systems}

Let $F$ be a non archimedean local field with finite residue class field $k_F$ of
characteristic $p$. Let $G={\rm GL}_{2}(F)$. Let $Z$ be the center of $G$. 

Let $X$ be the Bruhat-Tits tree of ${\rm 
PGL}_{2}/F$. Let $X^0$ denote its set of vertices, let $X^1$ denote its set of
edges; throughout, we identify an edge with its two-element set of
vertices. Let $d:X^0\times X^0\to{\mathbb Z}_{\ge0}$ be the
counting-edges-on-geodesics distance. By definition, an automorphism $g$ of
$X$ is a permutation $g$ of $X^0$ with $d(x,y)=d(gx,gy)$ for all $x,y\in X^0$;
clearly such a $g$ also induces also a
permutation $g$ of $X^1$.

We fix $\sigma=\{x_+,x_-\}\in X^1$.

For $x\in X^0$ let $K_x\subset G$ denote the maximal compact subgroup fixing
$x$, let $U_x^{(1)}\subset K_x$ denote its maximal normal
pro-$p$-subgroup. For $e\in{\mathbb N}$ let $U_x^{(e)}=\cap_{y}U_y^{(1)}$
where $y$ runs through all $y\in X^0$ with $d(x,y)\le e-1$. For
$\tau=\{x_1,x_2\}\in X^1$ let $U_{\tau}^{(e)}\subset G$ be the subgroup
generated by $U_{x_1}^{(e)}$ and $U_{x_2}^{(e)}$; this is again a
pro-$p$-group. For $e\ge0$ we put $$Z^{(e)}(x)=\{y\in X^0\,|\,d(x,y)\le
e\}\quad\quad\mbox{ for }x\in X^0.$$If $e\ge1$ then
$Z^{(e)}(x)=(X^0)^{U_{x}^{(e)}}$, the fixed point set of $U_{x}^{(e)}$ acting
on $X^0$. Next, for $e\ge0$ we put $$Z^{(e)}(\tau)=Z^{(e)}(x_1)\cup Z^{(e)}(x_2)\quad\quad\mbox{ for }\tau=\{x_1,x_2\}\in X^1.$$ We have $Z^{(e)}(\tau)=(X^0)^{U_{\tau}^{(e+1)}}$. \\

Let ${{\mathfrak o}}$ be a complete discrete valuation ring field with residue
class field $k$ of characteristic $p$. Let ${\rm Art}({\mathfrak o})$ denote
the category of Artinian local ${\mathfrak o}$-algebras with residue class
field $k$.\\

{\bf Definition:} Let $\Lambda\in{\rm Art}({\mathfrak o})$. A homological coefficient
system ${\mathcal F}$ in $\Lambda$-modules on $X$ is a collection of data as follows:

--- a $\Lambda$-module ${\mathcal F}(\tau)$ for each simplex $\tau$ 

--- a $\Lambda$-linear map $r^{\tau}_x:{\mathcal F}(\tau)\to {\mathcal F}(x)$ for each $x\in X^0$ and $\tau\in X^1$ with $x\in\tau$.

We obtain a $\Lambda$-linear map \begin{gather}\bigoplus_{\tau\in X^1}{\mathcal F}(\tau)\longrightarrow \bigoplus_{x\in X^0}{\mathcal F}(x)\label{difdef}\end{gather} sending $f\in{\mathcal F}(\tau)$ to $\sum_{x\in X^0\atop x\in\tau} r^{\tau}_x(f)$. The cokernel of the map (\ref{difdef}) is denoted by $H_0(X,{\mathcal F})$, its kernel is denoted by $H_1(X,{\mathcal F})$.

\begin{lem}\label{standardlahm} Let ${\mathcal F}$ be a homological coefficient system on $X$ for which the transition maps $r^{\tau}_x:{\mathcal F}(\tau)\to{\mathcal F}(x)$ are injective, for all $x\in\tau\in X^1$.

For any $y\in X^0$ the natural map ${\mathcal F}(y)\to H_0(X,{\mathcal F})$ is injective. In particular, if $H_0(X,{\mathcal F})=0$, then ${\mathcal F}=0$. 
\end{lem}

{\sc Proof:} Suppose that the map (\ref{difdef}) sends $c=(c_{\tau})_{\tau}\in
\oplus_{\tau\in X^1}{\mathcal F}(\tau)$ to the submodule ${\mathcal
  F}(y)$ of $\oplus_{x\in X^0}{\mathcal F}(x)$ (i.e. to an element with zero-component in all ${\mathcal
  F}(y')$ for all $y'\ne y$). We claim that $c=0$. Otherwise
there is some $\tau\in X^1$ with $c_{\tau}\ne0$ and some $x\in \tau$ such that
$d(y,x)$ is maximal (for all such $\tau$ and $x$). But then $x\ne y$ and the
injectivity of $r^{\tau}_x:{\mathcal F}(\tau)\to{\mathcal F}(x)$ shows
$c_{\tau}=0$, contradiction.\hfill$\Box$\\

Let $H$ be a group (or a monoid) acting on $X$ (through automorphisms of $X$). We say that the homological coefficient system ${\mathcal F}$ is $H$-equivariant if in addition we are given a $\Lambda$-linear map $g_{\tau}:{\mathcal F}(\tau)\to {\mathcal F}(g\tau)$ for each $g\in H$ and each ($0$- or
$1$-)simplex $\tau$, subject to the following conditions: 

(a) $g_{h\tau}\circ h_{\tau}=(gh)_{\tau}$ for each $g,h\in H$ and each simplex $\tau$ 

(b) $1_{\tau}={\rm id}_{{\mathcal F}(\tau)}$ for each simplex $\tau$ 

(c) $r^{g\tau}_{gx}\circ g_{x}=g_{\tau}\circ r_{x}^{\tau}$ for each $g\in H$ and each $x\in X^0$ and $\tau\in X^1$ with $x\in\tau$.\\

It is clear that if ${\mathcal F}$ is an $H$-equivariant homological coefficient system, then $H$ acts compatibly on the source and on the target of the map (\ref{difdef}), hence it acts on $H_0(X,{\mathcal F})$ and on $H_1(X,{\mathcal F})$. There is an obvious notion of a morphism ${\mathcal F}\to{\mathcal G}$ between $H$-equivariant homological coefficient systems: a collection of maps ${\mathcal F}(\tau)\to{\mathcal G}(\tau)$ for all simplices $\tau$, compatible with the restriction maps and the $H$-actions.\\
 
{\bf Definition:} For $\Lambda\in{\rm Art}({\mathfrak o})$ let ${\mathfrak
  C}^{(e)}(\Lambda)$ denote
the category of $G$-equivariant homological coefficient systems ${\mathcal F}$
in $\Lambda$-modules on $X$ satisfying the following conditions:\\(a) for any ($0$- or
$1$-)simplex $\tau$ the action of $U^{(e)}_{\tau}$ on ${\mathcal F}(\tau)$ is
trivial,\\(b) for any $z\in\eta\in X^1$ the transition map
$\tau_{z}^{\eta}:{\mathcal F}(\eta)\to{\mathcal F}(z)$ is injective, its image
is ${\mathcal F}(z)^{U^{(e)}_{\eta}}$, and this image ${\mathcal
  F}(z)^{U^{(e)}_{\eta}}$ generates ${\mathcal F}(z)$ as a
$K_z$-representation. [Thus, for $z\in X^0$, if $S=\{\eta\in X^1\,|\,z\in\eta\}$, then ${\mathcal F}(z)=\sum_{\eta\in S}{\rm im}(\tau_{z}^{\eta})$.]\\

Let $V$ be a $G$-representation on a $\Lambda$-module. We define a
coefficient system ${\mathcal F}_V^{(e)}$ on $X$ as
follows: $${\mathcal
  F}_V^{(e)}(\eta)=V^{{{U}}^{(e)}_{\eta}} \quad\quad\mbox{ and }\quad\quad{\mathcal
  F}_V^{(e)}(x)=\sum_{y\in X^0\atop \{x,y\}\in
  X^1}V^{{{U}}^{(e)}_{\{x,y\}}}=\sum_{y\in X^0\atop \{x,y\}\in
  X^1}{\mathcal
  F}^{(e)}_V(\{x,y\})$$ for $\eta\in
X^1$ and $x\in X^0$ (where in the definition of ${\mathcal
  F}_V^{(e)}(x)$ the sum is taken inside $V$).

\begin{lem} ${\mathcal F}_V^{(e)}$ belongs to ${\mathfrak
  C}^{(e)}(\Lambda)$. 
\end{lem}

{\sc Proof:} This is obvious.\hfill$\Box$\\

For $\Lambda\in{\rm Art}({\mathfrak o})$ let ${\mathfrak R}^{(e)}(\Lambda)$ denote the category of $G$-representations on $\Lambda$-modules which are
  generated by their ${{U}}^{(e)}_{\sigma}$-fixed vectors.

\begin{satz}\label{jalgcentral} Assume that $F={\mathbb Q}_p$. The assignments ${\mathcal F}\mapsto
  H_0(X,{\mathcal F})$ and $V\mapsto{\mathcal F}_V^{(1)}$ are an equivalence of categories between ${\mathfrak C}^{(1)}(\Lambda)$ and
  ${\mathfrak R}^{(1)}(\Lambda)$. For ${\mathcal F}\in {\mathfrak C}^{(1)}(\Lambda)$ and $V\in {\mathfrak R}^{(1)}(\Lambda)$ the natural maps \begin{gather}{\mathcal
  F}\longrightarrow {\mathcal F}_{H_0(X,{\mathcal F})}^{(1)}\quad\quad\mbox{ and }\quad\quad H_0(X,{\mathcal
    F}_V^{(1)})\longrightarrow V\end{gather}are isomorphisms.
\end{satz}

{\sc Proof:} \cite{jalg}.\hfill$\Box$\\

\section{Subgroups of the automorphism group of the tree}

$\quad\quad${\bf Definition:} For $e\in{\mathbb N}$ let $\widehat{G}^{(e)}$ denote the set of automorphisms
$g$ of $X$ with the property that for all $\mu\in X^1$ there is a $g'\in G$
such that the restrictions of $g$ and $g'$ to $Z^{(e)}(\mu)$ (viewed as maps
$Z^{(e)}(\mu)\to X^0$) coincide. $\widehat{G}^{(e)}$ is easily seen to be a
subgroup of ${\rm Aut}(X)$; we call it the group of locally algebraic
automorphisms of $X$ of level $e$. We have the chain of group inclusions$${\rm
  PGL}_2(F)\subset\ldots
\subset\widehat{G}^{(e+1)}\subset\widehat{G}^{(e)}\subset\ldots\subset
\widehat{G}^{(1)}\subset{\rm Aut}(X).$$

The following Theorem \ref{choucit} follows from the work of Choucroun \cite{choucrou}. This theorem (like the notations we need in order to formulate it) is not needed in the sequel, but of course it should be stated.

Fix a sequence of vertices $\ldots,x_{-2},x_{-1},x_0,x_1,x_2,\ldots$ forming a geodesic in $X$. Let $\widehat{ B}_+^{(e)}$ (resp. $\widehat{ B}_-^{(e)}$) denote the stabilizer in $\widehat{G}^{(e)}$ of the end of $X$ corresponding to $x_0,x_1,x_2,\ldots$ (resp. to $\ldots, x_{-2},x_{-1},x_0$). Let $\widehat{ N}_+^{(e)}$ denote the subgroup of $\widehat{ B}_+^{(e)}$ consisting of all $g\in \widehat{ B}_+^{(e)}$ with $g(x_n)=x_n$ for almost all $n\ge0$. Let $\widehat{K}^{(e)}$ denote the stabilizer in $\widehat{G}^{(e)}$ of the vertex $x_0$. Let $\widehat{I}^{(e)}$ denote the stabilizer in $\widehat{G}^{(e)}$ of the edge $\{x_0,x_1\}$. Let $\widehat{T}^{(e)}$ denote the pointwise stabilizer in $\widehat{G}^{(e)}$ of $\ldots,x_{-2},x_{-1},x_0,x_1,x_2,\ldots$.  

Let $s\in{\rm
  PGL}_2(F)\subset \widehat{G}^{(e)}$ be an element with $s(x_n)=x_{-n}$ for all $n\in{\mathbb Z}$, let $\varphi\in{\rm
  PGL}_2(F)\subset \widehat{G}^{(e)}$ be an element with $\varphi(x_n)=x_{n+1}$ for all $n\in{\mathbb Z}$.

\begin{satz}\label{choucit} (a) Cartan decomposition: $$\widehat{G}^{(e)}=\coprod_{m\in{\mathbb Z}_{\ge0}}\widehat{K}^{(e)}\varphi^m\widehat{K}^{(e)}.$$

 (b) Iwasawa decomposition: $$\widehat{G}^{(e)}=\widehat{K}^{(e)}\widehat{ B}_+^{(e)}.$$

(c) Bruhat decomposition: $$\widehat{G}^{(e)}=\widehat{N}_+^{(e)}s\widehat{ B}_+^{(e)}\coprod \widehat{ B}^{(e)}_+$$and moreover, $n_1s\widehat{ B}_+^{(e)}\cap n_1s\widehat{ B}_+^{(e)}\ne\emptyset$ if and only if $n_1\widehat{T}^{(e)}=n_2\widehat{T}^{(e)}$. 

(d) We have $$\widehat{ B}_+^{(e)}=\varphi^{\mathbb Z}\widehat{ N}_+^{(e)}.$$

(e) We have $$\widehat{ I}^{(e)}=(\widehat{ I}^{(e)}\cap\widehat{ B}_+^{(e)})\cdot(\widehat{ I}^{(e)}\cap\widehat{ B}_-^{(e)}),$$$$\widehat{T}^{(e)}=(\widehat{ I}^{(e)}\cap\widehat{ B}_+^{(e)})\cap(\widehat{ I}^{(e)}\cap\widehat{ B}_-^{(e)}).$$

\end{satz}

{\sc Proof:} This follows from \cite{choucrou} Theorem 1.5.2 as $\widehat{G}^{(e)}$ is closed in ${\rm Aut}(X)$ and contains the group ${\rm
  PGL}_2(F)$ which (in the terminology of \cite{choucrou}) acts weakly two transitive on $X$.\hfill$\Box$\\

{\bf Definition:} For a $1$-simplex $\eta$ let $\widehat{U}^{(e)}_{\eta}$ denote the set of
automorphisms $g$ of $X$ with the property that for all $\mu\in X^1$ there is
a $g'\in U_{\eta}^{(e)}$ such that the restrictions of $g$ and $g'$ to
$Z^{(e)}(\mu)$ (viewed as maps $Z^{(e)}(\mu)\to X^0$) coincide. Notice that
this implies: for all $x\in X^0$ there is
a $g'\in U_{\eta}^{(e)}$ such that the restrictions of $g$ and $g'$ to
$Z^{(e)}(x)$ coincide. \\

Clearly $\widehat{U}^{(e)}_{\eta}$ is a subgroup of $\widehat{G}^{(e)}$, and the
natural map $G\to\widehat{G}^{(e)}$ restricts to a map ${U}^{(e)}_{\eta}\to\widehat{U}^{(e)}_{\eta}$.\\

{\bf Remark:} We do {\it not} impose any topology on $\widehat{G}^{(e)}$. However, in the following we mimick the smooth representation theory of $p$-adic reductive groups by assigning to the subgroups $\widehat{U}^{(e)}_{\eta}$ of $\widehat{G}^{(e)}$ the role of open subgroups. On the other hand, on these individual subgroups $\widehat{U}^{(e)}_{\eta}$ we do consider their pro-$p$-topology (cf. Lemma \ref{nprop}).\\

\begin{lem}\label{nprop} ${\widehat{U}}^{(e)}_{\sigma}$ is a pro-$p$-group, for any $e\in{\mathbb N}$.
\end{lem}

{\sc Proof:} For any $m\ge0$, restriction induces a group homomorphism
${\widehat{U}}^{(e)}_{\sigma}\to{\rm Aut}(Z^{(m)}(\sigma))$ to the symmetric
group ${\rm Aut}(Z^{(m)}(\sigma))$ on the set $Z^{(m)}(\sigma)$, and we
have $${\widehat{U}}^{(e)}_{\sigma}={\rm lim}_{\leftarrow\atop m}{\rm
  im}[{\widehat{U}}^{(e)}_{\sigma}\longrightarrow {\rm
  Aut}(Z^{(m)}(\sigma))].$$Therefore we need to show that all the ${\rm
  im}[{\widehat{U}}^{(e)}_{\sigma}\to {\rm Aut}(Z^{(m)}(\sigma)]$ are finite
$p$-groups. We do this by induction on $m$. For $m<e$ the map
${\widehat{U}}^{(e)}_{\sigma}\to{\rm Aut}(Z^{(m)}(\sigma))$ is trivial. Now
let $m\ge e$. Let $g\in{\widehat{U}}^{(e)}_{\sigma}$ and $z\in
{Z^{(m-e)}(\sigma)}$. Then $g^{p^N}(z)=z$ for some $N\ge0$ by induction
hypothesis. As $g^{p^N}\in{\widehat{U}}^{(e)}_{\sigma}$ we find some $h\in
U_{\sigma}^{(e)}$ such that the restrictions of $h$ and $g^{p^N}$ to
$Z^{(e)}(z)$ coincide. Since $g^{p^N}(z)=z$ this is an equality $h=g^{p^N}$ in
${\rm Aut}(Z^{(e)}(z))$. Now $U_{\sigma}^{(e)}$ is a pro-$p$-group, therefore
we find some $M\ge0$ such that $h^{p^M}$ acts trivially on
$Z^{(m)}(\sigma)$. Therefore $g^{p^{N+M}}$ acts trivially on $Z^{(e)}(z)$. Thus, there is some $K\ge0$ such that $g^{p^{K}}$ acts
trivially on $Z^{(e)}(z)$ for any $z\in Z^{(m-e)}(\sigma)$. As $$Z^{(m)}(\sigma)\subset \bigcup_{z\in Z^{(m-e)}(\sigma)}Z^{(e)}(z)$$ we are done.\hfill$\Box$\\

{\bf Definition:} Let ${\mathfrak C}_0^{(e)}(\Lambda)$ denote the category of all
${\mathcal F}\in{\mathfrak C}^{(e)}(\Lambda)$ for which the $G$-action
factors through ${\rm PGL}_2(F)=G/Z$.\\

Let ${\mathcal
  F}\in{\mathfrak C}_0^{(e)}(\Lambda)$ and
$g\in{\widehat{G}}^{(e)}$. Given $\eta\in X^1$, choose a $g'\in G$ restricting
  to $g$ on $Z^{(e)}(\eta)$ and define $g_{\eta}:{\mathcal F}(\eta)\to{\mathcal
    F}(g\eta)$ to be the map $g'_{\eta}:{\mathcal F}(\eta)\to{\mathcal
    F}(g'\eta)={\mathcal F}(g\eta)$. This definition is independent of the
  choice of $g'$. Indeed, let also $g''\in G$ restrict to $g$ on
  $Z^{(e)}(\eta)$.  Then $g'^{-1}g''$ belongs to the pointwise stabilizer ${\rm
    Stab}_G(Z^{(e)}(\eta))$ of $Z^{(e)}(\eta)$ in $G$. If $x$ is one of the
  vertices of $\eta$, we have$${\rm
    Stab}_G(Z^{(e)}(\eta))\subset {\rm
    Stab}_G(Z^{(e)}(x))=U^{(e)}_xZ\subset U^{(e)}_{\eta}Z.$$Since
  $U^{(e)}_{\eta}Z$ acts trivially on ${\mathcal F}(\eta)$ by the
  definition of ${\mathfrak C}_0^{(e)}(\Lambda)$, we see that $g'^{-1}g''$ acts
  trivially on ${\mathcal F}(\eta)$. Therefore $g'$ and $g''=g'(g'^{-1}g'')$ define the same maps ${\mathcal F}(\eta)\to{\mathcal
    F}(g\eta)$. Similarly, given $x\in X^0$,
we choose a $g'\in G$ restricting to $g$ on $Z^{(e)}(x)$ and define
$g_x:{\mathcal F}(x)\to{\mathcal F}(gx)$ to be the map $g'_x:{\mathcal
  F}(x)\to{\mathcal F}(g'x)={\mathcal F}(gx)$; again this does not depend on the choice of $g'$.

\begin{lem}\label{actoncoef} The above definitions make ${\mathcal F}$ into a ${\widehat{G}}^{(e)}$-equivariant coefficient system on $X$. In particular, ${\widehat{G}}^{(e)}$ acts on $H_0(X,{\mathcal F})$.
\end{lem}

{\sc Proof:} This is clear.\hfill$\Box$\\

\begin{lem}\label{hatuinv} Let $\Lambda\in {\rm Art}({\mathfrak o})$. For any
  ${\mathcal F}\in{\mathfrak C}_0^{(e)}(\Lambda)$ the natural maps \begin{gather}{\rm im}(r_{x_{?}}^{\sigma})={\mathcal F}(x_{?})^{{{U}}^{(e)}_{\sigma}}\longrightarrow
  H_0({X},{\mathcal F})^{{\widehat{U}}^{(e)}_{\sigma}}\label{surjui}\end{gather}are bijective, for
  both $?=+$ and $?=-$.
\end{lem}

{\sc Proof:} The injectivity follows from Lemma \ref{standardlahm}. Now we prove surjectivity. Let the $0$-chain $c=(c_v)_{v\in X^0}$ represent the class $[c]$
in $H_0({X},{\mathcal F})^{{\widehat{U}}^{(e)}_{\sigma}}$. Let
$n(c)\in{\mathbb Z}_{\ge0}$ be minimal with ${\rm supp}(c)\subset
Z^{(n(c))}(\sigma)$. By induction on $n(c)$ we show that $[c]$ lies in the
image of the map (\ref{surjui}). If $n(c)=0$ the statement is clear: use the
injectivity of ${\mathcal F}(x_{?})\to H_0({X},{\mathcal F})$ (Lemma \ref{standardlahm}). Now let $n(c)>0$.

{\it Claim: We have $c_z\in{\mathcal F}(z)^{U^{(e)}_{\{z^-,z\}}}$ for all $z\in Z^{(n(c))}(\sigma)$, where $z^-\in Z^{(n(c)-1)}(\sigma)$ is such that $\{z^-,z\}\in X^1$.}

Given the claim, the defining properties of ${\mathfrak C}_0^{(e)}(\Lambda)$ allow us to pass from $c=(c_v)_{v\in X^0}$ to a homologous $0$-chain $c'=(c'_v)_{v\in X^0}$ with $n(c')=n(c)-1$. Then the induction hypothesis can be applied. 

To {\it prove the claim}, let $g\in {U^{(e)}_{\{z,z^-\}}}$. We find $g'\in
U_{z^-}^{(e)}$ with $gU_{z}^{(e)}=g'U_{z}^{(e)}$. Let $w\in Z^{(n(c)-2)}(\sigma)$ be such that $\{w,z^-\}\in X^1$. (If $n(c)=1$ and hence $z^-\in\sigma$, then take $w$ such that $\{w,z^-\}=\sigma$.) Removing $\{w,z^-\}$ from $X$ we are left with two disjoint closed full subhalftrees of $X$: the halftree $X_2$ rooted at $w$, and the half tree $X_1$ rooted at $z^{-}$. As $g'$ fixes $\{w,z^-\}$ pointwise, the action of $g'$ on $X$ respects $X_2$ and $X_1$. Let $\widehat{g}\in{\rm Aut}(X)$ be the unique element fixing $X_2$ pointwise and acting on $X_1$ like $g'$. It then follows by construction that in fact $\widehat{g}$ belongs to ${\widehat{U}}^{(e)}_{\sigma}$ and satisfies the following
properties: 

(i) For any $y\in X^0$ with $d(y,z^-)<d(y,z)$ we have $\widehat{g}(y)=y$, and
$\widehat{g}$ acts trivially on ${\mathcal F}(y)$.

(ii) We have $g'=\widehat{g}$ on ${\mathcal F}(z)$, and hence $g=\widehat{g}$
on ${\mathcal F}(z)$ (as
$U_{z}^{(e)}$ acts trivially on ${\mathcal F}(z)$).   

 As the support of $c$ is contained in the set of all $y\in X^0$ mentioned in
(i), together with $z$, we have$$\widehat{g}(c)-c=\widehat{g}(c_{z})-c_{z}={g}(c_{z})-c_{z}.$$On the
other hand, the class $[c]$ is ${\widehat{U}}^{(e)}_{\sigma}$-invariant, i.e. $\widehat{g}(c)-c$ maps to the
zero element in $H_0({X},{\mathcal F})$. Together, using Lemma
\ref{standardlahm} we see that ${g}(c_{z})=c_{z}$, as desired.\hfill$\Box$\\

{\bf Definition:} We define $\widehat{\mathfrak R}_0^{(e)}(\Lambda)$ to be the category of
representations of ${\widehat{G}}^{(e)}$ on $\Lambda$-modules which are generated by their ${\widehat{U}}^{(e)}_{\sigma}$-fixed vectors and which, when restricted to ${\widehat{U}}^{(e)}_{\sigma}$, are smooth ${\widehat{U}}^{(e)}_{\sigma}$-representations.\\

{\bf Definition:} Given $V\in \widehat{\mathfrak R}_0^{(e)}(\Lambda)$ we define a
coefficient system  $\widehat{\mathcal F}_V^{(e)}$ as
follows: $$\widehat{\mathcal
  F}_V^{(e)}(\eta)=V^{{\widehat{U}}^{(e)}_{\eta}} \quad\quad\mbox{ and }\quad\quad\widehat{\mathcal
  F}_V^{(e)}(x)=\sum_{y\in X^0\atop \{x,y\}\in
  X^1}V^{{\widehat{U}}^{(e)}_{\{x,y\}}}=\sum_{y\in X^0\atop \{x,y\}\in
  X^1}\widehat{\mathcal
  F}^{(e)}_V(\{x,y\})$$ for $\eta\in
X^1$ and $x\in X^0$ (where in the definition of $\widehat{\mathcal
  F}_V^{(e)}(x)$ the sum is taken inside $V$). The transition map
$r_x^{\eta}:\widehat{\mathcal F}_V^{(e)}(\eta)\to\widehat{\mathcal
  F}_V^{(e)}(x)$ for $x\in \eta$ is defined as follows: if $x$ lies in the
${\rm SL}_2(F)$-orbit ${\rm SL}_2(F)x_+$ of $x_+$, then $r_x^{\eta}$ is the
inclusion; if however $x\in {\rm SL}_2(F)x_-$ then $r_x^{\eta}$ is the {\it
  negative} of the inclusion. (Notice that $X^0={\rm SL}_2(F)x_+\coprod {\rm
  SL}_2(F)x_-$.) \\

The following result is an analogue of Theorem \ref{jalgcentral}.

\begin{satz}\label{centralequiv} (a) For ${\mathcal F}\in {\mathfrak
    C}_0^{(e)}(\Lambda)$ we have $H_0(X,{\mathcal F})\in \widehat{\mathfrak R}_0^{(e)}(\Lambda)$. For $V\in \widehat{\mathfrak R}_0^{(e)}(\Lambda)$ we have $\widehat{\mathcal F}_V^{(e)}\in {\mathfrak C}_0^{(e)}(\Lambda)$.

(b) These assignments ${\mathcal F}\mapsto
  H_0(X,{\mathcal F})$ and $V\mapsto\widehat{\mathcal F}_V^{(e)}$ are
  functorial in a natural way and form an adjoint pair: for ${\mathcal F}\in {\mathfrak C}_0^{(e)}(\Lambda)$ and $V\in\widehat{\mathfrak R}_0^{(e)}(\Lambda)$ we have a natural isomorphism\begin{gather}{\rm Hom}_{\widehat{\mathfrak R}_0^{(e)}(\Lambda)}(H_0(X,{\mathcal F}),V)\cong{\rm Hom}_{{\mathfrak C}_0^{(e)}(\Lambda)}({\mathcal F},\widehat{\mathcal F}_V^{(e)}).\label{adju}\end{gather}

(c) These functors are an
  equivalence of categories between ${\mathfrak C}_0^{(e)}(\Lambda)$ and
  $\widehat{\mathfrak R}_0^{(e)}(\Lambda)$. 

(d) For ${\mathcal F}\in {\mathfrak C}_0^{(e)}(\Lambda)$ and $V\in \widehat{\mathfrak R}_0^{(e)}(\Lambda)$ the natural maps \begin{gather}{\mathcal
  F}\longrightarrow \widehat{\mathcal F}_{H_0(X,{\mathcal F})}^{(e)}\label{fiso}\end{gather}\begin{gather}H_0(X,\widehat{\mathcal
    F}_V^{(e)})\longrightarrow V\label{viso}\end{gather}are isomorphisms.
\end{satz}

{\sc Proof:} (a) For
$\{x,z\}\in X^1$ we first claim that\begin{gather}({U}^{(e)}_{\{x,z\}},\bigcap_{y\in X^0\atop\{x,y\}\in
  X^1}\widehat{U}^{(e)}_{\{x,y\}})=\widehat{U}^{(e)}_{\{x,z\}}.\label{uinde}\end{gather}The
containment of the left hand side in $\widehat{U}^{(e)}_{\{x,z\}}$ is
obvious. Conversely, let $u\in\widehat{U}^{(e)}_{\{x,z\}}$. Multiplying by an
element in ${U}^{(e)}_{\{x,z\}}$ we may assume that $u$ fixes
$Z^{(e)}(x)$ pointwise. But then it is easy to see
that $u\in \widehat{U}^{(e)}_{\{x,y\}}$ for any $y\in X^0$ with $\{x,y\}\in
  X^1$. Now let $V\in \widehat{\mathfrak R}_0^{(e)}(\Lambda)$. We claim that\begin{gather}(\sum_{y\in X^0\atop\{x,y\}\in
  X^1}V^{\widehat{U}^{(e)}_{\{x,y\}}})^{{U}^{(e)}_{\{x,z\}}}=V^{\widehat{U}^{(e)}_{\{x,z\}}}.\label{vinde}\end{gather}The
containment of $V^{\widehat{U}^{(e)}_{\{x,z\}}}$ in the left hand side is
obvious. The reverse containment follows from formula (\ref{uinde}). Formula (\ref{vinde}) says $\widehat{\mathcal
  F}_V^{(e)}(x)^{{U}^{(e)}_{\{x,z\}}}=\widehat{\mathcal
  F}_V^{(e)}(\{x,z\})$. This shows that $\widehat{\mathcal F}_V^{(e)}\in
{\mathfrak C}_0^{(e)}(\Lambda)$.

Next, let ${\mathcal F}\in {\mathfrak C}_0^{(e)}(\Lambda)$. Since $H_0(X,{\mathcal F})$ is generated (as a ${\widehat{G}}^{(e)}$-representation) by ${\mathcal F}(\sigma)$, it is in particular generated by its ${\widehat{U}}^{(e)}_{\sigma}$-fixed vectors. Moreover, the ${\widehat{U}}^{(e)}_{\sigma}$-action on $H_0(X,{\mathcal F})$ is smooth: Indeed, given an element $v$ of $H_0(X,{\mathcal F})$, we pick a $0$-cycle $c\in C_0(X,{\mathcal F})$ representing $v$. We find $m\ge e$ such that $c$ is supported on $Z^{(m-e)}(\sigma)$. Thus $c$ and hence $v$ is fixed by the kernel of ${\widehat{U}}^{(e)}_{\sigma}\to{\rm Aut}(Z^{(m)}(\sigma))$; this kernel is of finite index in ${\widehat{U}}^{(e)}_{\sigma}$, hence open in ${\widehat{U}}^{(e)}_{\sigma}$ (cf. the proof of Lemma \ref{nprop}). We have shown that $H_0(X,{\mathcal F})\in \widehat{\mathfrak R}_0^{(e)}(\Lambda)$. 

For statement (b), the proof is the same as in \cite{jalg} Lemma 1.2. Statements (a), (b) and (d) together imply statement (c).

(d) We prove the bijectivity of the map (\ref{viso}). The composite $$V^{\widehat{U}_{\sigma}^{(e)}}=\widehat{\mathcal F}_V^{(e)}({\sigma})\stackrel{r_{x_{+}}^{\sigma}}{\longrightarrow}\widehat{\mathcal F}_V^{(e)}(x_{+})^{{\widehat{U}}^{(e)}_{\sigma}} {\longrightarrow} H_0({X},\widehat{\mathcal
  F}_V^{(e)})^{{\widehat{U}}^{(e)}_{\sigma}}{\longrightarrow}V$$is just the inclusion of $V^{\widehat{U}_{\sigma}^{(e)}}$ into $V$. The map $V^{\widehat{U}_{\sigma}^{(e)}}\to H_0({X},\widehat{\mathcal
  F}_V^{(e)})^{{\widehat{U}}^{(e)}_{\sigma}}$ is
surjective by Lemma \ref{hatuinv}, therefore the map $H_0({X},\widehat{\mathcal
  F}_V^{(e)})^{{\widehat{U}}^{(e)}_{\sigma}}\to V$ is injective. If the map
(\ref{viso}) was not injective, then, as (\ref{viso}) is
${\widehat{U}}^{(e)}_{\sigma}$-equivariant, the kernel of (\ref{viso}) would
contain a non-zero $\Lambda[{\widehat{U}}^{(e)}_{\sigma}]$-submodule. But
then, as $\Lambda$ is Artinian with residue field $k$, this kernel would also
contain a non-zero $k[{\widehat{U}}^{(e)}_{\sigma}]$-submodule. As
${\widehat{U}}^{(e)}_{\sigma}$ is a pro-$p$-group by Lemma \ref{nprop}, as it acts smoothly on $H_0(X,{\mathcal F})$ and as ${\rm char}(k)=p$, this kernel would therefore have a non-zero vector invariant under ${\widehat{U}}^{(e)}_{\sigma}$: contradiction to the injectivity of $H_0({X},\widehat{\mathcal
  F}_V^{(e)})^{{\widehat{U}}^{(e)}_{\sigma}}\to V$. Thus, the map (\ref{viso}) is injective. Its surjectivity is clear as $V$ is generated by its ${\widehat{U}}^{(e)}_{\sigma}$-invariant vectors.

Finally, it remains to prove that the map (\ref{fiso}) is an isomorphism. Fix
${\mathcal F}\in {\mathfrak C}_0^{(e)}(\Lambda)$. For the course of this proof
let us write ${\mathcal F}_{H}$ instead of $\widehat{\mathcal F}^{(e)}_{H_0(X,{\mathcal
    F})}$. For $\tau\in X^1$ let $x_{\tau,+}\in X^0$ denote the vertex of
$\tau$ belonging to the orbit ${\rm SL}_2(F)x_+$. As ${\mathcal F}$ has
injective transition maps the composition $${\mathcal
  F}(\tau)\stackrel{\tau_{x_{\tau,+}}^{\tau}}{\longrightarrow}{\mathcal
  F}(x_{\tau,+}){\longrightarrow}H_0(X,{\mathcal F})$$ is injective, see Lemma
\ref{standardlahm}. Clearly the map ${\mathcal F}_{H}(\tau)\stackrel{\tau_{x_{\tau,+}}^{\tau}}{\longrightarrow}{\mathcal
  F}_{H}(x_{\tau,+}){\longrightarrow}H_0(X,{\mathcal F})$ is
injective, too. Using these maps we regard ${\mathcal F}(\tau)$ and ${\mathcal F}_{H}(\tau)$ as
being contained in $H_0(X,{\mathcal F})$. Our hypotheses on ${\mathcal F}$
and the definition of ${\mathcal F}_{H}$ show that in this way we may
regard ${\mathcal F}$ as a sub coefficient system of ${\mathcal
  F}_{H}$. Namely, as $\widehat{U}^{(e)}$ acts trivially on
${\mathcal F}(\tau)$ for $\tau\in X^1$, we have the injective maps
$$\alpha_{\tau}:{\mathcal F}(\tau)\longrightarrow H_0(X,{\mathcal F})^{\widehat{U}^{(e)}_{\tau}}={\mathcal
  F}_{H}(\tau),$$ and as ${\mathcal F}(x)=\sum_{\tau\in
  X^1\atop x\in \tau}{\mathcal F}(\tau)$ the $\alpha_{\tau}$ also induce
injective maps $$\alpha_{x}:{\mathcal F}(x)=\sum_{\tau\in
  X^1\atop x\in \tau}{\mathcal F}(\tau)\longrightarrow \sum_{\tau\in
  X^1\atop x\in \tau}H_0(X,{\mathcal F})^{\widehat{U}^{(e)}_{\tau}}={\mathcal
  F}_{H}(x)$$ for $x\in X^0$. This morphism of coefficient systems ${\mathcal F}\to {\mathcal
  F}_{H}$ induces a map
$\alpha:H_0(X,{\mathcal F})\to H_0(X,{\mathcal F}_{H})$. On the other hand, let $\beta: H_0(X,{\mathcal
  F}_{H})\to H_0(X,{\mathcal F})$ be the natural morphism corresponding to the identity on ${\mathcal
  F}_{H}$ under the isomorphism (\ref{adju}).

{\it Claim: $\beta\circ\alpha$ is the identity on $H_0(X,{\mathcal F})$.} 

As $H_0(X,{\mathcal F})$ is generated by the images of the natural maps
$\iota_x:{\mathcal F}(x)\to H_0(X,{\mathcal F})$ for all $x\in X^0$, it
is enough to show $\beta\circ\alpha\circ\iota_x=\iota_x$ for all $x\in
X^0$. If $$\eta_x:{\mathcal
  F}_{H}(x)\longrightarrow H_0(X,{\mathcal
  F}_{H})$$ denotes the natural map, then we have
  $\alpha\circ\iota_x=\eta_x\circ \alpha_x$ by the definition of $\alpha$. Now
  by the definition of $\beta$ we have that $\beta\circ \eta_x$ is just the inclusion $${\mathcal
  F}_{H}(x)=\sum_{\tau\in
  X^1\atop x\in \tau}H_0(X,{\mathcal F})^{\widehat{U}^{(e)}_{\tau}}\to H_0(X,{\mathcal F})$$
for all $x\in X^0$. It follows that $\beta\circ\alpha\circ\iota_x=\beta\circ
\eta_x\circ \alpha_x$ is the inclusion of ${\mathcal F}(x)$ into $H_0(X,{\mathcal F})$, i.e. it is the
map $\iota_x$, as desired. The claim is proven. 

By the bijectivity of the map (\ref{viso}), applied to $V=H_0(X,{\mathcal F})$, the map $\beta$ is an isomorphism; hence $\alpha$ is an isomorphism, by the above claim. In particular, $H_0(X,{\mathcal F}_{H}/{\mathcal F})=0$. But it follows from our hypotheses on ${\mathcal F}$ that for all $x\in\eta\in X^1$ we have $${\mathcal F}_{H}(\eta)\cap{\mathcal F}(x)={\mathcal F}(\eta)$$inside ${\mathcal F}_{H}(x)$, i.e. $({\mathcal F}_{H}/{\mathcal F})(\eta)\to({\mathcal F}_{H}/{\mathcal F})(x)$ is injective, i.e. the quotient system ${\mathcal F}_{H}/{\mathcal F}$ has injective transition maps. These two facts together imply ${\mathcal F}_{H}/{\mathcal F}=0$, see Lemma \ref{standardlahm}. We get ${\mathcal F}_{H}={\mathcal F}$, i.e. the map (\ref{viso}) is an isomorphism.\hfill$\Box$\\

\section{Modules over the pro-$p$-Iwahori Hecke algebra}

 Let $\Lambda\in {\rm Art}({\mathfrak o})$ and
 $e\in{\mathbb N}$. Let ${\mathcal J}^{(e)}_{\Lambda}={\rm
   ind}_{U_{\sigma}^{(e)}}^{{G}}{\bf 1}_{\Lambda}$ denote the
 ${G}$-representation on the ${\Lambda}$-module of compactly supported
 ${\Lambda}$-valued functions on $U_{\sigma}^{(e)}\backslash{G}$. Let $${\mathcal H}_{\Lambda}(G,{U}^{(e)}_{\sigma})={\rm End}_G({\mathcal J}^{(e)}_{\Lambda}),$$the
 Hecke algebra of ${U}^{(e)}_{\sigma}\subset G$, with coefficients in
$\Lambda$. (For $e=1$ this is the pro-$p$-Iwahori Hecke algebra over $\Lambda$.) It is naturally isomorphic to
the $\Lambda$-algebra $\Lambda[{U}^{(e)}_{\sigma}\backslash G/{U}^{(e)}_{\sigma}]$ (in
which multiplication is given by convolution). This isomorphism
sends the coset
${U}^{(e)}_{\sigma}g{U}^{(e)}_{\sigma}$, for $g\in G$, to the endomorphism of ${\mathcal J}^{(e)}_{\Lambda}$
which sends the (compactly supported,
${U}^{(e)}_{\sigma}$-left invariant) function $f: G\to\Lambda$ to the
function $$G\longrightarrow\Lambda,\quad\quad h\mapsto\sum_{t\in{U}^{(e)}_{\sigma}\backslash
G}\chi_{{U}^{(e)}_{\sigma}g{U}^{(e)}_{\sigma}}(ht^{-1})f(t)$$(where
$\chi_{{U}^{(e)}_{\sigma}g{U}^{(e)}_{\sigma}}:G\to\Lambda$ denotes the
characteristic function of ${U}^{(e)}_{\sigma}g{U}^{(e)}_{\sigma}$). 

For any $G$-representation $V$ over $\Lambda$, the
$\Lambda$-algebra ${\mathcal H}_{\Lambda}(G,{U}^{(e)}_{\sigma})\cong
\Lambda[{U}^{(e)}_{\sigma}\backslash G/{U}^{(e)}_{\sigma}]$ naturally acts
from the right on the $\Lambda$-module $V^{{U}^{(e)}_{\sigma}}$ of
${U}^{(e)}_{\sigma}$-invariant vectors in $V$. The action of an arbitrary
coset ${U}^{(e)}_{\sigma}g{U}^{(e)}_{\sigma}$ on $v\in
V^{{U}^{(e)}_{\sigma}}$ is given by the following formula: if the collection $\{g_j\}_j$ in
$G$ is
such that
${U}^{(e)}_{\sigma}g{U}^{(e)}_{\sigma}=\coprod_j{U}^{(e)}_{\sigma}g_j$, then$$v\cdot{U}^{(e)}_{\sigma}g{U}^{(e)}_{\sigma}=\sum_jg_j^{-1}v.$$

\begin{pro}\label{heckestru} For any $V\in \widehat{\mathfrak R}_0^{(e)}(\Lambda)$ we have a
  natural right action of ${\mathcal H}_{\Lambda}(G,{U}^{(e)}_{\sigma})$ on
$V^{\widehat{U}^{(e)}_{\sigma}}$.
\end{pro}

{\sc Proof:} We may regard $V$ as a $G$-representation, therefore ${\mathcal
  H}_{\Lambda}(G,{U}^{(e)}_{\sigma})$ naturally acts on
$V^{{U}^{(e)}_{\sigma}}$. We claim that this action preserves the submodule
$V^{\widehat{U}^{(e)}_{\sigma}}$ of $V^{{U}^{(e)}_{\sigma}}$. 

By Theorem \ref{centralequiv} we have $H_0(X,{\mathcal
    F})\cong V$ for ${\mathcal
    F}=\widehat{\mathcal
    F}_V^{(e)}$. Let $v\in V^{\widehat{U}^{(e)}_{\sigma}}\cong H_0(X,{\mathcal
    F})^{\widehat{U}^{(e)}_{\sigma}}$, let $g\in G$ and $h\in\widehat{U}^{(e)}_{\sigma}$. We need to
      show \begin{gather}v\cdot{U}^{(e)}_{\sigma}g{U}^{(e)}_{\sigma}=h(v\cdot{U}^{(e)}_{\sigma}g{U}^{(e)}_{\sigma}).\label{uhatu}\end{gather}

First assume $g\sigma=\sigma$. By Lemma \ref{hatuinv} we may
represent $v$ by a $0$-cycle $c=(c_y)_{y\in X^0}\in C_0(X,{\mathcal F})$
supported on the vertices of $\sigma$. Choose $h'\in {U}^{(e)}_{\sigma}$ with
$h|_{Z^{(e)}(\sigma)}=h'|_{Z^{(e)}(\sigma)}$. As $c$ and hence $gc$ is supported on the vertices of $g\sigma=\sigma$ we see
$h(v\cdot{U}^{(e)}_{\sigma}g{U}^{(e)}_{\sigma})=h'(v\cdot{U}^{(e)}_{\sigma}g{U}^{(e)}_{\sigma})$,
therefore statement (\ref{uhatu}) follows from
$V^{\widehat{U}^{(e)}_{\sigma}}\subset V^{{U}^{(e)}_{\sigma}}$.

Now assume $g\sigma\ne\sigma$. Let $z\in\sigma$ be the vertex
of $\sigma$ such that, if $z'\in\sigma$ denotes the
other vertex, then $d(z,y)>d(z',y)$ for $y\in g\sigma$. Observing
$g{U}^{(e)}_{\sigma}g^{-1}={U}^{(e)}_{g\sigma}$ we easily see that\begin{gather}ZK_z\cap g{U}^{(e)}_{\sigma}g^{-1}\subset
  {U}^{(e)}_{\sigma}.\label{kcapgu}\end{gather}

For any collection $\{g_j\}_{j}$ in $g{U}^{(e)}_{\sigma}$ we deduce from this the
equivalence \begin{gather}{U}^{(e)}_{\sigma}g{U}^{(e)}_{\sigma}=\coprod_j{U}^{(e)}_{\sigma}g_j\,\,\,\mbox{ in }G\quad\quad\Leftrightarrow\quad\quad
{U}^{(e)}_{\sigma}g^{-1}z=\coprod_j\{g_j^{-1}z\}\,\,\,\mbox{ in
}X^0.\label{geohec}\end{gather}Indeed, as
$X^0\cong G/ZK_z$ the second statement in formula (\ref{geohec}) is equivalent to the statement
$ZK_zg{U}^{(e)}_{\sigma}=\coprod_jZK_zg_j$. But using formula (\ref{kcapgu}) this is straightforwardly
checked to be equivalent to the
first statement in formula (\ref{geohec}).

Choose a collection $\{g_j\}_{j}$ in $g{U}^{(e)}_{\sigma}$ satisfying the equivalent conditions of
formula (\ref{geohec}). Moreover, choose $h_j\in {U}^{(e)}_{\sigma}$ with
$h|_{Z^{(e)}(g_j^{-1}z)}=h_j|_{Z^{(e)}(g_j^{-1}z)}$, for any $j$. We then
find $${U}^{(e)}_{\sigma}g^{-1}z=\widehat{U}^{(e)}_{\sigma}g^{-1}z=h\widehat{U}^{(e)}_{\sigma}g^{-1}z=h{U}^{(e)}_{\sigma}g^{-1}z=\coprod_j\{hg_j^{-1}z\}=\coprod_j\{h_jg_j^{-1}z\}$$and
hence, in view of the equivalence (\ref{geohec}), also \begin{gather}{U}^{(e)}_{\sigma}g{U}^{(e)}_{\sigma}=\coprod_j{U}^{(e)}_{\sigma}g_jh_j^{-1}.\label{halt}\end{gather} 

By Lemma \ref{hatuinv} we may
represent $v$ by a $0$-cycle $c=(c_y)_{y\in X^0}\in C_0(X,{\mathcal F})$
supported on $z$. Then $v\cdot{U}^{(e)}_{\sigma}g{U}^{(e)}_{\sigma}$ is
represented by the $0$-cycle $c^g=(c^g_y)_{y\in X^0}$ with
$c^g_{g_j^{-1}z}=g_j^{-1}c_{z}$ for all $j$, and with $c^g_y=0$ for all other
vertices $y$. Therefore $h(v\cdot{U}^{(e)}_{\sigma}g{U}^{(e)}_{\sigma})$ is
represented by the $0$-cycle $h(c^{g})=((hc^{g})_y)_{y\in X^0}$ with
$(hc^{g})_{h_jg_j^{-1}z}=h_jg_j^{-1}c_{z}$ for all $j$, and $(hc^{g})_y=0$ for all other
vertices $y$. But by formula (\ref{halt}) this is again a representative for $v\cdot{U}^{(e)}_{\sigma}g{U}^{(e)}_{\sigma}$.\hfill$\Box$\\

Let $W$ be a $ZK_{x_+}$-representation on a free $\Lambda$-module such that $U^{(1)}_{x_+}$
acts trivially on $W$. We assume that for any $z\in X^0$ with $\{x,z\}\in X^1$
the natural map\begin{gather}W\longrightarrow\bigoplus_{y\in
    X^0-\{z\}\atop\{x,y\}\in
    X^1}W_{U^{(1)}_{\{x,y\}}}\label{injinco}\end{gather}is injective. Here $W_{U^{(1)}_{\{x,y\}}}$ denotes the module of coinvariants for the action of $U^{(1)}_{\{x,y\}}$ on $W$. Consider
the compact induction ${\rm ind}_{ZK_{x_+}}^GW$, the $\Lambda$-module of all
locally constant functions $f:G\to W$, with compact support modulo $Z$, which
satisfy $f(gk)=k^{-1}f(g)$ for all $k\in ZK_{x_+}$. The group $G$ acts by left
translation on ${\rm ind}_{ZK_{x_+}}^GW$. For $x\in X^0$ we choose $g\in G$
with $x=gx_+$ and put $${\mathcal G}(x)={\mathcal G}_W(x)=\{f\in{\rm
  ind}_{ZK_{x_+}}^GW\,|\,{\rm supp}(f)\subset gZK_{x_+}\}$$(this does not
depend on the choice of $g$, as $X^0\cong G/ZK_{x_+}$). We
have \begin{gather}{\rm ind}_{ZK_{x_+}}^GW=\bigoplus_{x\in X^0}{\mathcal
    G}(x).\label{heckeop}\end{gather}Suppose that we are given a
$G$-equivariant endomorphism $T$ of ${\rm ind}_{ZK_{x_+}}^GW$ which, as an
endomorphism of $\oplus_{x\in X^0}{\mathcal G}(x)$ (via the identification
(\ref{heckeop})), has the following structure: for any $x\in X^0$, the
restriction $T|_{{\mathcal G}(x)}$ factors as \begin{gather}T|_{{\mathcal G}(x)}:{\mathcal
  G}(x)\stackrel{T_x}{\longrightarrow}  \bigoplus_{y\in X^0\atop \{x,y\}\in X^1}{\mathcal
  G}(y)\stackrel{\iota}{\hookrightarrow}\bigoplus_{y\in X^0}{\mathcal
  G}(y)\label{tstru}\end{gather}where $\iota$ is the natural inclusion, and where $T_x$ is the sum of
maps\begin{gather}T_{x,y}:{\mathcal G}(x){\longrightarrow} {\mathcal
  G}(x)_{U^{(1)}_{\{x,y\}}}\stackrel{t_{x,y}}{\longrightarrow}{\mathcal
  G}(y)^{U^{(1)}_{\{x,y\}}}\longrightarrow{\mathcal G}(y)\label{ttstru}\end{gather}for all $y\in X^0$
with $\{x,y\}\in X^1$, where the first arrow, resp. the last arrow, is the
canonical projection onto the coinvariants, resp. the canonical inclusion of
the invariants, and where $t_{x,y}$ is an isomorphism of $\Lambda$-modules. 

\begin{lem}\label{stanhec} Let $\lambda\in \Lambda$. If for $b=\sum_{x\in X^0}b_x\in \oplus_{x\in X^0}{\mathcal
    G}(x)$ the support of $$(T-\lambda)(b)=\sum_xT(b)_x-\lambda b_x\in \oplus_{x\in X^0}{\mathcal
    G}(x)$$ is contained in the two vertices of some $\eta\in X^1$, (i.e. $T(b)_x=\lambda b_x$ for
  $x\notin \eta$), then $b=0$.
\end{lem}

{\sc Proof:} This is parallel to Lemma \ref{standardlahm}. Assume
$b\ne0$. Choose $x\in X^0$ with $b_x\ne0$ such that $d(\eta,x)={\rm min}\{d(y,x)\,|\,y\in\eta\}$ is maximal (for
all such $x$). Then $T(b)_{x'}=T_{x,x'}(b_x)$ for all $x'\in X^0$ with
$\{x,x'\}\in X^1$ and $d(x',\eta)=d(x,\eta)+1$. The definition of the $T_{x,x'}$ and the injectivity of the map
(\ref{injinco}) implies, by translation, the injectivity of the map$$\sum_{x'}T_{x,x'}:{\mathcal
  G}(x)\longrightarrow\bigoplus_{x'}{\mathcal
  G}(x')$$where the sums are taken over all $x'\in X^0$ with $\{x,x'\}\in X^1$ and $d(x',\eta)=d(x,\eta)+1$. Therefore $T(b)_{x'}=0$ for all such $x'$ implies $b_x=0$:
contradiction. \hfill$\Box$\\

\begin{pro}\label{packet} Let $\lambda\in\Lambda$. Assume that $Z$ acts trivially on $W$.

(a) $\widehat{G}^{(1)}$ naturally acts on ${\rm ind}_{ZK_{x_+}}^GW$, and $T$
is $\widehat{G}^{(1)}$-equivariant. In particular, $\widehat{G}^{(1)}$, and
hence its subgroup $\widehat{U}^{(1)}_{\sigma}$, act on $({\rm ind}_{ZK_{x_+}}^GW)/{}(T-\lambda)$. 

(b) The image of ${\mathcal G}(x_+)^{{U}^{(1)}_{\sigma}}\oplus{\mathcal G}(x_-)^{{U}^{(1)}_{\sigma}}$ in $({\rm ind}_{ZK_{x_+}}^GW)/{}(T-\lambda)$ is contained in the submodule of $\widehat{U}^{(1)}_{\sigma}$-invariants, and the natural map \begin{gather}{\mathcal G}(x_+)^{{U}^{(1)}_{\sigma}}\bigoplus{\mathcal G}(x_-)^{{U}^{(1)}_{\sigma}}\longrightarrow (\frac{{\rm ind}_{ZK_{x_+}}^GW}{{}(T-\lambda)})^{\widehat{U}^{(1)}_{\sigma}}\label{breuilkey}\end{gather}is bijective. 
\end{pro}

{\sc Proof:} (Here and below we write $({\rm ind}_{ZK_{x_+}}^GW)/{}(T-\lambda)$ instead of $({\rm ind}_{ZK_{x_+}}^GW)/{\rm
  im}(T-\lambda)$). 

(a) The action of $\widehat{G}^{(1)}$ on ${\rm ind}_{ZK_{x_+}}^GW$ is defined in the same way as the action of $\widehat{G}^{(e)}$ on ${\mathcal
  F}\in{\mathfrak C}_0^{(e)}(\Lambda)$, cf. Lemma \ref{actoncoef}. The $\widehat{G}^{(1)}$-equivariance of $T$ is then immediate from its local nature.

(b) The injectivity even of ${\mathcal G}(x_+)\oplus{\mathcal G}(x_-)\to ({\rm ind}_{ZK_{x_+}}^GW)/{}(T-\lambda)$ follows from Lemma \ref{stanhec}. The following proof of surjectivity of the map (\ref{breuilkey}) is similar to that of Lemma \ref{hatuinv}. Let $a=\{a_x\}_{x\in X^0}\in\bigoplus_{x\in X^0}{\mathcal G}(x)={\rm ind}_{ZK_{x_+}}^GW$ be such that its class $[a]$ in $({\rm ind}_{ZK_{x_+}}^GW)/{}(T-\lambda)$ is $\widehat{U}^{(1)}_{\sigma}$-invariant. Let $$n(a)={\rm max}\{d(\sigma,x)\,|\,x\in X^0, a_x\ne0\}$$where we put $d(\sigma,x)={\rm min}\{d(x_+,x), d(x_-,x)\}$.

{\it Step 1:} Let $x\in X^0$ with $d(\sigma,x)=n(a)-1$. We claim that for any
element $x'$ of the set $${\mathcal O}_x=\{x'\in X^0\,\,|\,\,\{x,x'\}\in X^1,
d(\sigma,x')=n(a)\}$$we have $a_{x'}\in {\mathcal
  G}(x')^{U^{(1)}_{\{x,x'\}}}$. 

To see this, let $u\in U^{(1)}_{\{x,x'\}}$. We find some $u'\in U^{(1)}_x$
with $u\equiv u'$ modulo $U^{(1)}_{x'}$. Next, we find some $\widehat{u}\in
\widehat{U}^{(1)}_{\sigma}$ such that 

(i) for all $y\in X^0$ with $d(y,x)<d(y,x')$ we have $\widehat{u}(y)=y$, and
$\widehat{u}$ acts trivially on ${\mathcal
  G}(y)$, and

(ii) we have $u'=\widehat{u}$ on ${\mathcal
  G}(x')$, and hence $u=\widehat{u}$ on ${\mathcal
  G}(x')$ (as $U^{(1)}_{x'}$ acts trivially on ${\mathcal
  G}(x')$).

As the support of $a$ is contained in the set consisting of $x'$ and all $y\in X^0$ mentioned in
(i), we
have $$\widehat{u}(a)-a=\widehat{u}(a_{x'})-a_{x'}={u}(a_{x'})-a_{x'}.$$On the
other hand, the class $[a]$ is $\widehat{U}^{(1)}_{\sigma}$-invariant,
i.e. $\widehat{u}(a)-a\in{}(T-\lambda)$. Together with Lemma \ref{stanhec} we
obtain ${u}(a_{x'})=a_{x'}$, as desired.

{\it Step 2:} Let $x$ and ${\mathcal O}_x$ be as in step 1. We claim that
$\sum_{x'\in{\mathcal O}_x}a_{x'}$ lies in the image of the natural map
\begin{gather}{\mathcal G}(x)\longrightarrow\bigoplus_{x'\in{\mathcal O}_x}{\mathcal G}(x')_{U^{(1)}_{\{x,x'\}}}.\label{barliv}\end{gather}

By step 1, we may view $a_{x'}$, for any $x'\in{\mathcal O}_x$, by means of the isomorphism $t_{x,x'}: {\mathcal
  G}(x)_{U^{(1)}_{\{x,x'\}}}\cong{\mathcal
  G}(x')^{U^{(1)}_{\{x,x'\}}}$ as an element in ${\mathcal
  G}(x)_{U^{(1)}_{\{x,x'\}}}$. Let $x''\in X^0$ be the unique vertex with
$\{x,x''\}\in X^1$ and $d(\sigma,x'')\le d(\sigma,x)$. Then $U^{(1)}_{x''}$
acts transitively on the set ${\mathcal O}_x$. Given an element $u$ of
$U^{(1)}_{x''}$, we find an element $\widehat{u}$ of
$\widehat{U}^{(1)}_{\sigma}$ whose action on $\{y\in X^0\,|\,d(\sigma,y)\le
n(a)\}$ fixes (pointwise) all elements not contained in ${\mathcal O}_x$, and
acts on ${\mathcal O}_x$ like $u$. Applying Lemma
\ref{stanhec} again (and using the
$\widehat{U}^{(1)}_{\sigma}$-invariance of the class $[a]$) we therefore see that $\sum_{x'\in{\mathcal O}_x}a_{x'}$
is a $U^{(1)}_{x''}$-invariant element of the $U^{(1)}_{x''}$-representation
$\oplus_{x'\in{\mathcal O}_x}{\mathcal G}(x)_{U^{(1)}_{\{x,x'\}}}$ (on which
$U^{(1)}_{x''}$ acts by permuting the summands transitively). Therefore
$\sum_{x'\in{\mathcal O}_x}a_{x'}$ lies in the image of the map
(\ref{barliv}).

{\it Step 3:} By what we saw in Step 2 we may pass to another representative
(modulo ${}(T-\lambda)$) of $[a]$ which has zero contribution at all $x'\in
{\mathcal O}_x$. Doing this for all $x\in X^0$ with $d(\sigma,x)=n(a)-1$ we
obtain a representative $\widetilde{a}$ of $[a]$ with $n(\widetilde{a})=n(a)-1$. Proceeding by
induction on $n(a)$ we finally see that any
$\widehat{U}^{(1)}_{\sigma}$-invariant class in $({\rm
  ind}_{ZK_{x_+}}^GW)/{}(T-\lambda)$ is represented by some $a$ with $n(a)=0$,
i.e. by some $a\in{\mathcal
  G}(x_+)\oplus{\mathcal G}(x_-)$. As in Step 1 we then see that this
representative in fact belongs to ${\mathcal G}(x_+)^{{U}^{(1)}_{\sigma}}\oplus{\mathcal G}(x_-)^{{U}^{(1)}_{\sigma}}$. \hfill$\Box$\\

Assume that $k$ is algebraically closed. Let $W$ be an irreducible $k[K_{x_+}]$-module on which the center of $K_{x_+}$ acts trivially. We regard $W$ as a $k[ZK_{x_+}]$-module by letting $Z$ act trivially.

\begin{satz}\label{gleichgut} (a) The map (\ref{injinco}) is injective, and the Hecke operator $T$ on ${\rm ind}_{ZK_{x_+}}^GW$ constructed in \cite{bali} has the structure described above, i.e. is given by maps (\ref{tstru}), (\ref{ttstru}).

(b) For $\lambda\in k$ the ${\mathcal H}_{k}(G,{U}^{(1)}_{\sigma})$-module $(\frac{{\rm ind}_{ZK_{x_+}}^GW}{{}(T-\lambda)})^{\widehat{U}^{(1)}_{\sigma}}$ has $k$-dimension $2$. If it is irreducible, then $\frac{{\rm ind}_{ZK_{x_+}}^GW}{{}(T-\lambda)}$ is an irreducible $\widehat{G}^{(1)}$-representation. 

\end{satz}

{\sc Proof:} (a) See \cite{bali}. 

(b) Let $0\ne V\subset\frac{{\rm ind}_{ZK_{x_+}}^GW}{{}(T-\lambda)}$ be a
$\widehat{G}^{(1)}$-sub representation. The same argument as in the proof of Theorem \ref{centralequiv} shows that the action of $\widehat{U}^{(1)}_{\sigma}$ on ${\rm ind}_{ZK_{x_+}}^GW$, and hence on $\frac{{\rm ind}_{ZK_{x_+}}^GW}{{}(T-\lambda)}$ and on $V$, is smooth. Therefore, since $\widehat{U}^{(1)}_{\sigma}$ is a pro-$p$-group by Lemma \ref{nprop}, we have $V^{\widehat{U}^{(1)}_{\sigma}}\ne0$. Replacing $V$ by its
$\widehat{G}^{(1)}$-sub representation generated by
$V^{\widehat{U}^{(1)}_{\sigma}}$ we may assume that $V$ belongs to
$\widehat{\mathfrak R}_0^{(1)}(k)$, as does $\frac{{\rm
    ind}_{ZK_{x_+}}^GW}{{}(T-\lambda)}$. Proposition \ref{heckestru} says that
${\mathcal H}_{k}(G,{U}^{(1)}_{\sigma})$ acts on $(\frac{{\rm
    ind}_{ZK_{x_+}}^GW}{{}(T-\lambda)})^{\widehat{U}^{(1)}_{\sigma}}$,
respecting $V^{\widehat{U}^{(1)}_{\sigma}}$. Therefore we have
$V^{\widehat{U}^{(1)}_{\sigma}}=(\frac{{\rm
    ind}_{ZK_{x_+}}^GW}{{}(T-\lambda)})^{\widehat{U}^{(1)}_{\sigma}}$ if
$(\frac{{\rm ind}_{ZK_{x_+}}^GW}{{}(T-\lambda)})^{\widehat{U}^{(1)}_{\sigma}}$
is irreducible. This implies $V=\frac{{\rm ind}_{ZK_{x_+}}^GW}{{}(T-\lambda)}$
as $\frac{{\rm ind}_{ZK_{x_+}}^GW}{{}(T-\lambda)}$ is generated by $(\frac{{\rm ind}_{ZK_{x_+}}^GW}{{}(T-\lambda)})^{\widehat{U}^{(1)}_{\sigma}}$.

By Proposition \ref{packet} we may identify $(\frac{{\rm ind}_{ZK_{x_+}}^GW}{{}(T-\lambda)})^{\widehat{U}^{(1)}_{\sigma}}$ with ${\mathcal G}(x_+)^{{U}^{(1)}_{\sigma}}\oplus{\mathcal G}(x_-)^{{U}^{(1)}_{\sigma}}$. It is well known that ${\mathcal G}(x_+)^{{U}^{(1)}_{\sigma}}=W^{{U}^{(1)}_{\sigma}}$ and hence also ${\mathcal G}(x_-)^{{U}^{(1)}_{\sigma}}$ is one-dimensional, therefore $(\frac{{\rm ind}_{ZK_{x_+}}^GW}{{}(T-\lambda)})^{\widehat{U}^{(1)}_{\sigma}}$ is two dimensional.\hfill$\Box$\\

{\bf Definition:} A $\widehat{G}^{(1)}$-representation on a $k$-vector space $V$ is called admissible if $V^{\widehat{U}^{(1)}_{\sigma}}$ is a finite dimensional $k$-vector space.

\begin{pro}\label{quot} For any admissible irreducible
  $\widehat{G}^{(1)}$-representation on a $k$-vector space $V$ there exists
  some $\lambda\in k$ and an irreducible $k[K_{x_+}]$-module $W$ on which the
  center of $K_{x_+}$ acts trivially, and a surjective homomorphism of
  $\widehat{G}^{(1)}$-representations \begin{gather}\frac{{\rm
      ind}_{ZK_{x_+}}^GW}{{}(T-\lambda)}\longrightarrow V.\label{gastri}\end{gather}
\end{pro}

{\sc Proof:} As $V$ is admissible, $\sum_{\eta\in X^1\atop x_+\in\eta}V^{\widehat{U}^{(1)}_{\eta}}$ is a finite dimensional $k[K_{x_+}]$-module. Therefore it admits an irreducible $k[K_{x_+}]$-module-submodule $W$; the center of $K_{x_+}$ acts trivially. By the definition of $\widehat{\mathcal
    F}_V^{(1)}$ we may regard $W$ also as a $k[ZK_{x_+}]$-submodule of $\widehat{\mathcal
    F}_V^{(1)}(x_+)$, and this induces a natural homomorphism of $\widehat{G}^{(1)}$-representations ${\rm ind}_{ZK_{x_+}}^GW\to \oplus_{x\in X^0}\widehat{\mathcal
    F}_V^{(1)}(x)$. Consider the composition of $\widehat{G}^{(1)}$-representations $$\gamma:{\rm ind}_{ZK_{x_+}}^GW\longrightarrow\oplus_{x\in X^0}\widehat{\mathcal
    F}_V^{(1)}(x)\longrightarrow H_0(X,\widehat{\mathcal
    F}_V^{(1)})\longrightarrow V$$(we know but do not need here that $H_0(X,\widehat{\mathcal
    F}_V^{(1)})\to V$ is an isomorphism (Theorem \ref{centralequiv})). The ring ${\rm End}({\rm ind}_{ZK_{x_+}}^GW)$ acts
  (by precomposing) from the right on the finite dimensional $k$-vector space
  $${\rm Hom}_{G}({\rm ind}_{ZK_{x_+}}^GW,V)={\rm Hom}_{ZK_{x_+}}(W,V)={\rm
    Hom}_{ZK_{x_+}}(W,\sum_{\eta\in X^1\atop
    x_+\in\eta}V^{\widehat{U}^{(1)}_{\eta}}).$$Consider the action of $T\in
  {\rm End}({\rm ind}_{ZK_{x_+}}^GW)$ (in fact we have  $k[T]={\rm End}({\rm
    ind}_{ZK_{x_+}}^GW)$, see \cite{bali}) on the  $k[T]$-submodule of ${\rm
    Hom}_{G}({\rm ind}_{ZK_{x_+}}^GW,V)$ generated by $\gamma$. As $k$ is algebraically closed, this action of $T$ has as a
  non-zero eigen vector, providing a non-zero $G$-equivariant homomorphism
  (\ref{gastri}) for some $\lambda\in k$. But since $\gamma$ and $T$ are
  $\widehat{G}^{(1)}$-equivariant, also this map (\ref{gastri}) is
  $\widehat{G}^{(1)}$-equivariant. Its surjectivity follows from the irreducibility of $V$. \hfill$\Box$\\

The simple ${\mathcal H}_{k}(G,{U}^{(1)}_{\sigma})$-right modules have been
classified by Marie France Vign\'{e}ras in \cite{vigcomp}. In particular, in 
\cite{vigcomp} par. 3.2 it has been defined what it means for a simple
${\mathcal H}_{k}(G,{U}^{(1)}_{\sigma})$-right module to be {\it supersingular}. We
say that a ${\mathcal H}_{k}(G,{U}^{(1)}_{\sigma})$-right module has {\it trivial
central character} if for any $g\in Z$ the double coset
${U}^{(1)}_{\sigma}g{U}^{(1)}_{\sigma}\in {k}[{U}^{(1)}_{\sigma}\backslash G/{U}^{(1)}_{\sigma}]\cong{\mathcal
  H}_{k}(G,{U}^{(1)}_{\sigma})$ acts as the identity. (This definition
makes no reference to the center of the ring ${\mathcal
  H}_{k}(G,{U}^{(1)}_{\sigma})$; for the (easy) description of that center see
however \cite{vigcomp}.)

\begin{kor}\label{heckmain} The functor $V\mapsto V^{\widehat{U}^{(1)}_{\sigma}}$ from
  $\widehat{G}^{(1)}$-representations to ${\mathcal
    H}_{k}(G,{U}^{(1)}_{\sigma})$-right modules induces a bijection between the
  isomorphism classes of admissible irreducible
  $\widehat{G}^{(1)}$-representations and the isomorphism classes of simple
  ${\mathcal H}_{k}(G,{U}^{(1)}_{\sigma})$-right modules with trivial central character.
\end{kor}

{\sc Proof:} First, in \cite{vigcomp} Proposition 4.9 it is shown that the
functor $V\mapsto V^{{U}^{(1)}_{\sigma}}$ from ${G}$-representations
to ${\mathcal H}_{k}(G,{U}^{(1)}_{\sigma})$-right modules induces a bijection
between the isomorphism classes of the irreducible subquotients of principal
(series) representations of $G$, and the isomorphism classes of simple
${\mathcal H}_{k}(G,{U}^{(1)}_{\sigma})$-right modules which are {\it not}
supersingular. Therefore, to see that the non-supersingular simple
${\mathcal H}_{k}(G,{U}^{(1)}_{\sigma})$-right modules with trivial central character also ly in the image of our functor under discussion, it is enough to see that the irreducible subquotients of principal
(series) representations of $G$ on which the center $Z$ of $G$
acts trivially are in fact $\widehat{G}^{(1)}$-representations. In \cite{bali} it is shown that the irreducible principal
(series) representations of $G$ are
isomorphic to $\frac{{\rm
      ind}_{ZK_{x_+}}^GW}{{}(T-\lambda)}$ for some $\lambda\in k^{\times}$ and some irreducible $k[K_{x_+}]$-module $W$. Proposition \ref{packet} shows that $\frac{{\rm
      ind}_{ZK_{x_+}}^GW}{{}(T-\lambda)}$ is indeed a $\widehat{G}^{(1)}$-representation (if $Z$ acts trivially). Moreover, in \cite{bali} it is also
  shown that the principal
(series) representations of $G$ which are not irreducible are necessarily of
length two, and the subquotients are the following: a twist of the
one-dimensional $G$-representation, and a twist of the
Steinberg representation. Clearly the trivial (one-dimensional)
$G$-representation is a $\widehat{G}^{(1)}$-representation, too. But it is also well known that the Steinberg
representation of $G$ is isomorphic to $H_0(X,{\mathcal F})$ for some
${\mathcal F}\in{\mathfrak C}_0^{(1)}(k)$, hence is a $\widehat{G}^{(1)}$-representation. (For example, for $x\in X^0$ the
$K_x$-representation ${\mathcal F}(x)$ is the $|k_F|$-dimensional Steinberg
representation of the reductive quotient of $K_x$.)

Next, we claim that any supersingular simple ${\mathcal
  H}_{k}(G,{U}^{(1)}_{\sigma})$-right module $M$ with trivial central character
is of the form $V^{\widehat{U}^{(1)}_{\sigma}}$ for some admissible irreducible
  $\widehat{G}^{(1)}$-representation $V$. Indeed, by \cite{padi} Proposition 2.18 (which
  is a quotation from \cite{vigcomp}), the constructions in \cite{padi}
  Definition 6.2, Lemma 6.3 assign to $M=M_{\gamma}$ a coefficient system
  ${\mathcal F}={\mathcal V}_{\gamma}$
  on $X$ (notations $M_{\gamma}$ and ${\mathcal V}_{\gamma}$ from loc.cit). It
  is immediate from these constructions that ${\mathcal F}\in{\mathfrak
    C}^{(1)}(k)$, and even ${\mathcal
    F}\in{\mathfrak C}_0^{(1)}(k)$ if $M$ has trivial central
  character. Moreover, in that case, Lemma
  \ref{hatuinv} and \cite{padi} Lemma 6.4 show that $M\cong H_0(X,{\mathcal
    F})^{\widehat{U}^{(1)}_{\sigma}}$ as ${\mathcal
  H}_{k}(G,{U}^{(1)}_{\sigma})$-right modules.

Finally, it remains to show that for any admissible irreducible
  $\widehat{G}^{(1)}$-representation $V$ the ${\mathcal
    H}_{k}(G,{U}^{(1)}_{\sigma})$-right module $V^{\widehat{U}^{(1)}_{\sigma}}$ is
  simple. We use Proposition \ref{quot}, i.e. we choose a surjection
  (\ref{gastri}). If $\lambda\ne0$ then $\frac{{\rm
      ind}_{ZK_{x_+}}^GW}{{}(T-\lambda)}$ is a principal (series) representation
  (by \cite{bali}), and our previous disussion applies. If $\lambda=0$ then
  the disussion in \cite{padi} (in particular Lemma 6.1, Lemma 6.3) shows that
  the left hand side of the map (\ref{breuilkey}) is a supersingular simple ${\mathcal
  H}_{k}(G,{U}^{(1)}_{\sigma})$-right module. By Theorem \ref{gleichgut} it
is isomorphic to $V^{\widehat{U}^{(1)}_{\sigma}}$. \hfill$\Box$\\

\section{A variation of $(\varphi,\Gamma)$-modules}

\subsection{A functor ${\bf D}$}

We fix $e\in{\mathbb N}$ and a uniformizer $p_F\in{\mathcal O}_F$. \\

{\bf Definition:} We define the submonoids$$P^{(e)}_{\star}=\{\left( \begin{array}{cc}a&b\\0&1\end{array}\right)\,|\,a\in
{\mathcal O}_F-\{0\}, b\in p_F^{e-1}{\mathcal O}_F\}\quad\mbox{ and }\quad
T_{\star}=\{\left( \begin{array}{cc}a&0\\0&1\end{array}\right)\,|\,a\in
{\mathcal O}_F-\{0\}\}$$and the
subgroups$$N^{(e)}_0=\{\left( \begin{array}{cc}1&b\\0&1\end{array}\right)\,|\,
b\in p_F^{e-1}{\mathcal O}_F
\}\quad\mbox{ and }\quad N=\{\left( \begin{array}{cc}1&b\\0&1\end{array}\right)\,|\,
b\in F\}$$and$${T}=\{\left( \begin{array}{cc}a&0\\0&1\end{array}\right)\,|\,a\in
F^{\times}\}$$of $G$. We have $N_0^{(e)}\rtimes T_{\star}=P^{(e)}_{\star}$. We write $t=\left( \begin{array}{cc}
    p_F&0\\0&1\end{array}\right)\in T_{\star}$. We then have$${T}=\bigcup_{m\in{\mathbb N}}t^{-m}T_{\star}.$$ We assume that our fixed central vertices $x_+$ and $x_-$ are those with
$K_{x_+}={\rm GL}_2({\mathcal O}_F)$ and $x_-=t^{-1}x_+$. Recall that
$\sigma=\{x_+,x_-\}$. In the following, the end $\alpha_0=t^{-\infty}x_+$ of
$X$ will play a distinguished role. For any $x\in X^0$ there is a unique infinite (without backtracking) path $[\alpha_0,z]$ in $X$ starting at $z$ and passing through almost all $t^{-m}x_+$ with $m\ge0$. For $e'\ge0$ and $x\in X^0$ we define the subset$$Z_+^{(e')}(x)=\{z\in X^0\,|\,d(x,z)\le e' \mbox{ and }x\in
[\alpha_0,z]\}$$of $X^0$. We define $\widehat{N}_0^{(e)}$ to be the
subgroup of ${\rm Aut}(X)$ consisting of all $g\in{\rm Aut}(X)$ with the
property that for all $x\in X^0$ there is
a $g'\in {N}^{(e)}_0$ such that the restrictions of $g$ and $g'$ to
$Z_+^{(e)}(x)$ coincide. For $\mu\in X^1$ we let $x_+^{\mu}\in X^0$ be that vertex of $\mu$ such that the other vertex of $\mu$
belongs to $[\alpha_0,x_+^{\mu}]$. For example, $x_+^{\sigma}=x_+$.

For $x\in X^0$ we let $N_x^{(e)}$ be the subgroup of $N$ consisting of all
$g\in N$ which fix $Z_+^{(e-1)}(x)$ pointwise. For example, $N_{x_+}^{(e)}=N_0^{(e)}$. 

\begin{lem}\label{nauchprop} (a) $\widehat{N}_0^{(e)}$ is a (non-abelian)
  pro-$p$-group.

(b) For any $g\in \widehat{N}_0^{(e)}$ and any $\mu\in X^1$ there is
a $g'\in {N}_0^{(e)}$ such that the restrictions of $g$ and $g'$ to the
subset $Z^{(e)}(\mu)$ of $X^0$ coincide.

(c) For any $\mu\in X^1$ we have $N_{x_+^{\mu}}^{(e)}=U_{\mu}^{(e)}\cap N$.
\end{lem}

{\sc Proof:} (a) The same as for Lemma \ref{nprop}. Alternatively, one may use (b)
to see that $\widehat{N}_0^{(e)}$ is a subgroup of $\widehat{U}^{(e)}_{\sigma}$ and then just
quote the result of Lemma \ref{nprop}. Or one may just use the description of
$\widehat{N}_0^{(e)}$ given in Lemma \ref{propexp} below.

(b) Apply the defining condition for
$\widehat{N}_0^{(e)}$ to all $z\in [\alpha_0,x_+^{\mu}]$ with
$d(z,x_+^{\mu})\le e$.

(c) This is clear.\hfill$\Box$\\

We define the subgroup (cf. the proof of Lemma \ref{svz} below) $$\widehat{N}^{(e)}=\bigcup_{m\in{\mathbb
    N}}t^{-m}\widehat{N}^{(e)}_0t^{m}$$of ${\rm Aut}(X)$. We read
${T}$ as a subgroup of ${\rm Aut}(X)$ and define $\widehat{P}^{(e)}$ to be the
subgroup of ${\rm Aut}(X)$ generated by ${T}$ and
$\widehat{N}^{(e)}$. We define $\widehat{P}_{\star}^{(e)}$ to be the submonoid
of $\widehat{P}^{(e)}$ generated by $\widehat{N}_0^{(e)}$ and $T_{\star}$. 

\begin{lem}\label{svz} (a) We
  have $\widehat{N}_0^{(e)}\rtimes T_{\star}=\widehat{P}^{(e)}_{\star}$ and $\widehat{N}^{(e)}\rtimes {T}=\widehat{P}^{(e)}$. 

(b) For any $m\in{\mathbb N}$ and $s\in \widehat{P}^{(e)}_{\star}$ we have $s^m\widehat{N}^{(e)}_0s^{-m}\subset
\widehat{N}^{(e)}_0$. Moreover, $$\bigcap_{m\in{\mathbb N}}t^m\widehat{N}^{(e)}_0t^{-m}=1.$$

\end{lem}

{\sc Proof:} (a) This follows from $N_0^{(e)}\rtimes T_{\star}=P^{(e)}_{\star}$.

(b) For the first claim it is enough to check $s^m\widehat{N}^{(e)}_0s^{-m}\subset
\widehat{N}^{(e)}_0$ for $s\in T_{\star}$. For such $s$ we have $s^m{N}^{(e)}_0s^{-m}\subset
{N}_0^{(e)}$ for all $m\in{\mathbb N}$, and the claim follows. The second one
follows from $\cap_{m\in{\mathbb N}}t^m{N}^{(e)}_0t^{-m}=1$.\hfill$\Box$\\

{\bf Remark:} ${T}N$ (resp. $N$, resp. $N_0^{(e)}$) is naturally a subgroup of
$\widehat{P}^{(e)}$ (resp. of $\widehat{N}^{(e)}$, resp. of $\widehat{N}^{(e)}_0$), and similarly, $P^{(e)}_{\star}$ is naturally a submonoid of $\widehat{P}^{(e)}_{\star}$. One
might ask for retractions in the reverse direction. Let ${\mathfrak
  E}_P$ denote the set of ends of $X$ except for the end
$\alpha_0$. Let $\overline{\alpha}_0=t^{\infty}x_+=t^{\infty}x_-$ denote the end of $X$ which together with $\alpha_0$ spans the $T$-stable apartment in $X$. Then $\widehat{P}^{(e)}$ acts on ${\mathfrak
  E}_P$ and sending $n\in N$ to
$n\overline{\alpha}_0$ is a
bijection between $N$ ($\cong F$) and ${\mathfrak
  E}_P$. For any $g\in \widehat{P}^{(e)}$ and $\alpha\in {\mathfrak
  E}_P$ there is a unique $g_{\alpha}\in {T}N$ such that $g$ and
$g_{\alpha}$ coincide on all the vertices on the apartment $[\alpha_0,\alpha]$ from $\alpha_0$ to
$\alpha$. This defines for each $\alpha\in {\mathfrak
  E}_P$ a map \begin{gather}\widehat{P}^{(e)}\longrightarrow
{T}N,\quad\quad g\mapsto g_{\alpha}\label{retra}\end{gather}such
that for any $g\in
{T}N\subset \widehat{P}^{(e)}$ we have $g=g_{\alpha}$. Similarly,
the maps (\ref{retra}) are retractions $\widehat{N}^{(e)}\to N$ and $\widehat{P}^{(e)}_{\star}\to P^{(e)}_{\star}$ and $\widehat{N}^{(e)}_0\to
N_0^{(e)}$. Beware that they are {\it not} group homomorphisms. Rather, we
have$$(g\cdot f)_{\alpha}=g_{f(\alpha)}\cdot f_{\alpha}$$for $f,g\in
\widehat{P}^{(e)}$ and $\alpha\in {\mathfrak
  E}_P$.\\ 

We define the subset$$X_{+}^0={P}^{(1)}_{\star}.x_+,\quad\quad\quad\quad\mbox{
  resp. }\quad\quad\quad\quad X_{+}^1={P}^{(1)}_{\star}.\sigma$$ of $X^0$, resp. of
$X^1$. If $X_{+}$ denotes the maximal connected full closed subcomplex of $X$
such that $x_+$ belongs to its set of vertices but $x_-$ does not, then
$X_{+}^0$ is precisely the set of vertices of $X_+$, and $X_{+}^1-\{\sigma\}$
is its set of edges. To clarify: $\sigma=\{x_+,x_-\}$ belongs to the set which
we denote by $X_{+}^1$, but it is {\it not} an edge of the {\it closed}
subcomplex $X_+$ of $X$.

For any $n\in N/N_0^{(1)}$ the subset $nX^0_+$ of $X^0$ is stable under the
action of $\widehat{N}_0^{(e)}$. Restriction to this subset defines a quotient
$\widehat{N}_{0,n}^{(e)}$ of $\widehat{N}_0^{(e)}$
contained in the group of permutations of the set $nX^0_+$.

\begin{lem}\label{decpropexp} The natural homomorphism\begin{gather}\widehat{N}_0^{(e)}\longrightarrow\prod_{n\in
    N/N_0^{(1)}}\widehat{N}_{0,n}^{(e)}\label{prodn}\end{gather}is bijective. The factors $\widehat{
  N}_{0,n}^{(e)}$ are pairwise isomorphic.

\end{lem}

{\sc Proof:} The union $\cup_{n\in
    N/N_0^{(1)}}nX^0_+$ is disjoint in $X^0$ and stable under
  $\widehat{N}_0^{(e)}$, and the restriction map from $\widehat{N}_0^{(e)}$ to
  the group of permutations of $\cup_{n\in
    N/N_0^{(1)}}nX^0_+$ is injective, because any element of $\widehat{N}_0^{(e)}$
  acts trivially on $X^0-\cup_{z\in
    N/N_0^{(1)}}nX^0_+$. \hfill$\Box$\\     

{\bf Remark:} The index of $t^m\widehat{N}^{(e)}_0t^{-m}$ in $\widehat{N}^{(e)}_0$ is not finite for $m\in{\mathbb N}$. Lemma \ref{svz} says
that, except for this failure, the axiomatic of section 3 in
\cite{scviza} is satisfied. The following statement might be viewed as a natural replacement of the finiteness axiom in
\cite{scviza}: we have$$[\widehat{N}^{(e)}_{0,1}:\widehat{N}^{(e)}_{0,1}\cap
t^m\widehat{N}^{(e)}_{0}t^{-m}]=\prod_{j=0}^{m-1}q^{e+j}$$for any $m\ge1$,
where $q=|k_F|$; here we view the factor $\widehat{N}^{(e)}_{0,1}$ of $\widehat{N}^{(e)}_{0}$ as a
subgroup in $\widehat{N}^{(e)}_{0}$. Following the lines of \cite{scviza} one
may thus define a corresponding notion of \'{e}taleness. However, although the
point in our entire discussion is that one should consider
$\widehat{N}^{(e)}_0$-actions instead of just ${N}^{(e)}_0$-actions, the good
notion of \'{e}taleness seems to be just the usual one, as given in formulae
(\ref{classet1}) and (\ref{classet2}) below (i.e. based on the cosets
$t^i{N}^{(1)}_{0}t^{-i}$ in ${N}^{(1)}_{0}$).\\ 

Let $\Lambda\in {\rm Art}({\mathfrak o})$. We are interested in $NT$-equivariant homological coefficient systems in $\Lambda$-modules ${\mathcal F}$ on $X$
satisfying the following hypotheses:

(Hyp 1) for any $x\in\eta\in X^1$ the transition map
$\tau_{x}^{\eta}:{\mathcal F}(\eta)\to{\mathcal F}(x)$ is injective.

(Hyp 2) for any $\mu\in
X^1$, if $S_{\mu}=\{\eta\in X^1\,|\,\mu\cap \eta=\{x_+^{\mu}\}\}$, then \begin{gather}{\mathcal
  F}(x_+^{\mu})=\sum_{\eta\in S_{\mu}}{\rm im}(\mu_{x_+^{\mu}}^{\eta})\label{humphr}.\end{gather}

(Hyp 3) for any $\mu\in
X^1$, the image of $\tau_{x_+^{\mu}}^{\mu}:{\mathcal F}(\mu)\to{\mathcal F}(x_+^{\mu})$ is
  ${\mathcal F}(x_+^{\mu})^{N_{x_+^{\mu}}^{(e)}}$.\\

Before proceeding, we first provide examples for such ${\mathcal F}$. 

\begin{lem}\label{altpetdis} Let $N(k_F)$, $\overline{N}(k_F)$ be the unipotent radicals of opposite Borel subgroups in ${\rm GL}_2(k_F)$, let $W$ be a $k[{\rm GL}_2(k_F)]$-module which is generated by its subspace $W^{N(k_F)}$ of $N(k_F)$-invariants. Then $W$ is generated by $W^{N(k_F)}$ even as a $k[\overline{N}(k_F)]$-module.
\end{lem}

{\sc Proof:} If $W$ is finitely generated, then, by Nakayama's Lemma, applied to the local ring $k[\overline{N}(k_F)]$, to show that $W^{N(k_F)}$ generates $W$ as a $k[\overline{N}(k_F)]$-module it is enough to show this after reduction modulo the maximal ideal, i.e. it is enough to show that the composition $W^{N(k_F)}\to W\to W_{\overline{N}(k_F)}$ is surjective. If $W$ is an irreducible $k[{\rm GL}_2(k_F)]$-module this is easily verified. Next suppose that $W$ is a princial series representation, i.e. parabolically induced from a character of a Borel subgroup. In this case $W$ is isomorphic, as a $k[\overline{N}(k_F)]$-module, with a direct sum of the trivial one-dimensional $k[\overline{N}(k_F)]$-module and a copy of $k[\overline{N}(k_F)]$. In particular, $W_{\overline{N}(k_F)}$ is two dimensional. Similarly we see that $W^{N(k_F)}$ is two dimensional. The dual of $W$ is a again a principal series representation $W'$, in particular $(W')^{\overline{N}(k_F)}$ maps surjectively onto the space of $\overline{N}(k_F)$-invariants of the unique irreducible quotient of $W'$. Dually we obtain that the space of $\overline{N}(k_F)$-coinvariants of the unique irreducible subobject of $W$ maps injectively into the space of $\overline{N}(k_F)$-coinvariants of $W$. Thus applying what we said above to the unique irreducible subobject and to the unique irreducible quotient of $W$ we see that the composition $W^{N(k_F)}\to W\to W_{\overline{N}(k_F)}$ is surjective. 

Now suppose we are given a general $W$. Then, as $W$ is generated by $W^{N(k_F)}$, we see that $W$ is isomorphic to a quotient of a direct sum of
principal series representations of ${\rm GL}_2(k_F)$. Their property of being generated by their $N(k_F)$-invariants clearly passes to direct sums and quotients, hence to $W$.\hfill$\Box$\\

\begin{lem}\label{vorslev} Any ${\mathcal F}\in{\mathfrak
  C}^{(1)}(k)$ satisfies hypotheses (Hyp 1), (Hyp 2) and (Hyp 3) (with $e=1$).
\end{lem}

{\sc Proof:} (We do not ask that ${\mathcal F}\in{\mathfrak
  C}_0^{(1)}(k)$: the center $Z$ of $G$ may act nontrivially on
${\mathcal F}$.) The validity of hypotheses (Hyp 1) and (Hyp 3) is clear, only hypothesis (Hyp 2) requires a proof. The
$K_{x_+^{\mu}}$-action on ${\mathcal F}(x_+^{\mu})$ factors through the
quotient of $K_{x_+^{\mu}}$ isomorphic to ${\rm GL}_2(k_F)$. Thus,
abstractly ${\mathcal F}(x_+^{\mu})$ is isomorphic to a representation of ${\rm GL}_2(k_F)$ generated by its invariants under
the unipotent radical $N(k_F)$ of a Borel subgroup in ${\rm GL}_2(k_F)$. But then Lemma \ref{altpetdis} tells us that ${\mathcal F}(x_+^{\mu})$ is generated by its invariants under $N(k_F)$ even as a representation of the unipotent radical of an opposite Borel subgroup. Given the other properties of ${\mathcal F}\in{\mathfrak
  C}_0^{(1)}(k)$, this is precisely the meaning of formula (\ref{humphr}).\hfill$\Box$\\

We now fix an arbitrary ${\mathcal F}$ over an arbitrary $\Lambda\in {\rm Art}({\mathfrak o})$ satisfying the above conditions (a), (b), (c).

Let $g\in
\widehat{P}^{(e)}$. Given $\eta\in X^1$, choose a $g'\in NT$ restricting
  to $g$ on $Z_+^{(e)}(x_+^{\eta})$ and define $g:{\mathcal F}(\eta)\to{\mathcal
    F}(g\eta)$ to be the map $g'_{\eta}:{\mathcal F}(\eta)\to{\mathcal
    F}(g'\eta)={\mathcal F}(g\eta)$. Similarly, given $x\in X^0$, choose a $g'\in NT$ restricting
  to $g$ on $Z_+^{(e)}(x)$ and define $g:{\mathcal F}(x)\to{\mathcal
    F}(gx)$ to be the map $g'_{\eta}:{\mathcal F}(x)\to{\mathcal
    F}(g'x)={\mathcal F}(gx)$.

\begin{lem}\label{nactoncoef} The above action of
  $\widehat{P}^{(e)}$ on ${\mathcal F}$ is well defined and makes
  ${\mathcal F}$ into a $\widehat{P}^{(e)}$-equivariant coefficient
  system on $X$. In particular, $\widehat{P}^{(e)}$ naturally acts on
  $H_0(X,{\mathcal F})$.
\end{lem}

{\sc Proof:} The same as for Lemma \ref{actoncoef}.\hfill$\Box$\\

Since $X_+$ is closed in $X$ we have a natural map $H_0(X_+,{\mathcal F})\to
H_0(X,{\mathcal F})$. As in Lemma \ref{standardlahm} this map is seen to be
injective. Moreover, as $X_+$ is a $\widehat{P}^{(e)}_{\star}$-stable subcomplex of
$X$ we have $\widehat{P}^{(e)}_{\star}$ acting also on the submodule $H_0(X_+,{\mathcal F})$ of $H_0(X,{\mathcal F})$. The
corresponding $\widehat{N}_{0}^{(e)}$-action on $H_0(X_+,{\mathcal F})$ in fact factors through the quotient $\widehat{N}_{0,1}^{(e)}$ of $\widehat{N}_{0}^{(e)}$.

\begin{lem}\label{nauchinv} The action of $\widehat{N}_{0,1}^{(e)}$ on $H_0(X,{\mathcal F})$ is smooth.

(b) The natural map ${\mathcal F}(\sigma)\cong {\mathcal F}(x_+)^{{N}_0^{(e)}}\to H_0(X_+,{\mathcal F})^{\widehat{N}^{(e)}_{0,1}}$ is bijective.
\end{lem}

{\sc Proof:} (a) The same as for Theorem \ref{centralequiv} (a).

(b) The same as for Lemma \ref{hatuinv}; notice that $(N_0^{(e)},U_{x_+}^{(e)})=U_{\sigma}^{(e)}$.\hfill$\Box$\\

We now assume that ${\mathfrak o}$ is the ring of integers in a finite
extension field $L$ of ${\mathbb Q}_p$ with residue class field $k=k_L$ and
uniformizer $p_L\in{\mathfrak o}$. We assume that $\Lambda\in{\rm
  Art}({\mathfrak o})$ is a {\it quotient} of ${\mathfrak o}$. For a profinite
group $H$ we write$${\mathcal
  L}(H)={\mathfrak o}[[H]]$$for its completed group
ring (Iwasawa algebra) over ${\mathfrak o}$. For example, the natural
maps$${N}_0^{(e)}\stackrel{\iota}{\longrightarrow}\widehat{N}_0^{(e)}\stackrel{{\rm
  pr}}{\longrightarrow}\widehat{N}_{0,1}^{(e)}$$--- both ${\iota}$ and ${\rm
pr}\circ {\iota}$ are injective ---  induce corresponding morphisms
of Iwasawa algebras ${\mathcal
  L}({N}_0^{(e)})\to{\mathcal
  L}(\widehat{N}_0^{(e)})\to{\mathcal
  L}(\widehat{N}_{0,1}^{(e)})$. Using Lemma \ref{nauchinv} we obtain that the Pontrjagin dual$$D({\mathcal F})={\rm Hom}^{\rm ct}_{{\mathfrak o}}(H_0(X_+,{\mathcal
  F}),L/{\mathfrak o})$$of $H_0(X_+,{\mathcal
  F})$ is a module over ${\mathcal
  L}(\widehat{N}_{0,1}^{(e)})$. 

\begin{lem}\label{topnak} If the $k$-vector space ${\mathcal F}(x_+)^{{N}^{(e)}_0, p_L=0}$ has dimension $n<\infty$, then the ${\mathcal
  L}(\widehat{N}_{0,1}^{(e)})$-module $D({\mathcal F})$ can be
  generated by $n$ elements. 
\end{lem} 

{\sc Proof:} ${\mathcal
  L}(\widehat{N}_{0,1}^{(e)})$ is a local ring (as $\widehat{N}_{0,1}^{(e)}$ is a pro-$p$-group), its augmentation (left-)ideal $I$ is maximal, and $I^n\to0$. We may therefore apply the topological
Nakayama Lemma (see \cite{baho}) to the profinite ${\mathcal
  L}(\widehat{N}_{0,1}^{(e)})$ (left-)module $D({\mathcal F})$: it says that $D({\mathcal F})$ can be generated over ${\mathcal
  L}(\widehat{N}_{0,1}^{(e)})$ by $n$ elements ($n<\infty$) if its reduction
modulo $I$ can be generated by $n$ elements. By duality, this is the
case if $H_0(X_+,{\mathcal
  F})^{\widehat{N}_{0,1}^{(e)}, p_L=0}$ has dimension $\le n$. Now apply Lemma \ref{nauchinv}.\hfill$\Box$\\

We also define the ${\mathcal
  L}(\widehat{N}_{0,1}^{(e)})$-module$$D'({\mathcal F})={\rm Hom}^{\rm ct}_{{\mathfrak o}}(\frac{H_0(X_+,{\mathcal
  F})}{{\mathcal F}(\sigma)},L/{\mathfrak o})\cong {\rm Hom}^{\rm ct}_{{\mathfrak o}}(\frac{H_0(X,{\mathcal
  F})}{H_0(X_-,{\mathcal
  F})},L/{\mathfrak o}).$$Here we view ${\mathcal F}(\sigma)$ as embedded into $H_0(X_+,{\mathcal F})$ via ${\mathcal
  F}(\sigma)\to {\mathcal F}(x_+)$, and $H_0(X_-,{\mathcal
  F})$ is the $0$-homology group of ${\mathcal
  F}$ on $X_-$, the full {\it closed} subcomplex of $X$ with vertex set
$X_-^0=X^0-X_+^0$. (Thus, $\sigma$ belongs neither to the complex $X_+$ nor to the complex $X_-$.) 

\begin{lem} $(\widehat{P}_{\star}^{(e)})^{-1}$ naturally acts on $D({\mathcal F})$,
  and $\widehat{P}_{\star}^{(e)}$ naturally acts on $D'({\mathcal F})$.
\end{lem}

{\sc Proof:} As the monoid $\widehat{P}^{(e)}_{\star}$ acts on $H_0(X_+,{\mathcal
  F})$, the inverse monoid $(\widehat{P}^{(e)}_{\star})^{-1}$ acts
on its dual $D({\mathcal F})$. On the other hand, $\widehat{P}^{(e)}_{\star}$ itself acts on $D'({\mathcal F})$. Indeed, the group $\widehat{N}^{(e)}_{0,1}$ acts on $H_0(X_+,{\mathcal
  F})$, respecting the submodule ${\mathcal F}(\sigma)$, hence it acts on $D'({\mathcal F})={\rm Hom}^{\rm ct}_{{\mathfrak o}}(\frac{H_0(X_+,{\mathcal
  F})}{{\mathcal F}(\sigma)},L/{\mathfrak o})$. The monoid ${T}_{\star}^{-1}$ acts on $H_0(X,{\mathcal
  F})$, respecting the submodule $H_0(X_-,{\mathcal
  F})$, hence ${T}_{\star}$ acts on $D'({\mathcal F})\cong {\rm Hom}^{\rm ct}_{{\mathfrak o}}(\frac{H_0(X,{\mathcal
  F})}{H_0(X_-,{\mathcal
  F})},L/{\mathfrak o})$.\hfill$\Box$\\

Let ${\mathcal O}^{(e)}_{{\mathcal E}}$ denote the $p$-adic completion of the
localization of ${\mathcal L}(N^{(e)}_0)$ with respect to the complement of
$p_L{\mathcal L}(N_0^{(e)})$. \\

{\bf Definition:} An \'{e}tale
$((\widehat{P}_{\star}^{(e)})^{-1},\widehat{P}^{(e)}_{\star})$-pair over
${\mathcal L}(\widehat{N}^{(e)}_{0,1})$ is an injective morphism of ${\mathcal
  L}(\widehat{N}^{(e)}_{0,1})$-modules $\partial:D'\to D$, together with an action of the monoid
$(\widehat{P}_{\star}^{(e)})^{-1}$ on $D$ and of the monoid
$\widehat{P}_{\star}^{(e)}$ on $D'$, either of them extending the given action of the group $\widehat{N}^{(e)}_{0}$ (through its
quotient $\widehat{N}^{(e)}_{0,1}$). The following properties are required:

(a) For $m\in{\mathbb N}$ let $\psi_{t^m}$ denote the
endomorphism $t^{-m}$ of $D$, and let $\varphi_{t^m}$ denote the
endomorphism $t^{m}$ of $D'$. For any $m\in{\mathbb N}$ we
have\begin{gather}\psi_{t^m}\partial\varphi_{t^m}=\partial.\label{classet1}\end{gather}

(b) The cokernel of $\partial$ is finite. As ${\mathcal O}^{(e)}_{{\mathcal E}}$ is flat over ${\mathcal
  L}({N}_0^{(e)})$ (and as $\partial$ is in particular ${\mathcal
  L}({N}^{(e)}_0)$-linear)  it follows that $\partial$ induces an isomorphism$${\mathcal O}^{(e)}_{{\mathcal E}}\otimes_{{\mathcal
  L}({N}^{(e)}_0)}D'\stackrel{\cong}{\longrightarrow} {\mathcal O}^{(e)}_{{\mathcal E}}\otimes_{{\mathcal
  L}({N}^{(e)}_0)}D.$$We identify its source and its target
into a single object ${\bf D}$.

(c) Denote again by $\varphi_{t^m}$ and $\psi_{t^m}$ the endomorphisms of
${\bf D}$ obtain from those on $D$ and $D'$ by base extension ${\mathcal O}^{(e)}_{{\mathcal E}}\otimes_{{\mathcal
  L}({N}_0^{(e)})}(.)$. Then we have\begin{gather}\sum_{n\in N_0^{(1)}/t^{m}N_0^{(1)}t^{-m}}n\varphi_{t^m}\psi_{t^m}n^{-1}={\rm id}_{{\bf
    D}}.\label{classet2}\end{gather}

\begin{satz}\label{fet} Suppose that the $k$-vector space ${\mathcal F}(x_+)^{{N}^{(e)}_0,
    p_L=0}$ has dimension $n<\infty$. Then $(D'({\mathcal
  F})\to D({\mathcal F}))$ is an \'{e}tale
$((\widehat{P}_{\star}^{(e)})^{-1},\widehat{P}^{(e)}_{\star})$-pair over
${\mathcal L}(\widehat{N}^{(e)}_{0,1})$. Moreover, $D({\mathcal F})$ can be
  generated as an
  ${\mathcal L}(\widehat{N}^{(e)}_{0,1})$-module by $n$ elements.
\end{satz}

{\sc Proof:} The last statement was shown in
  Lemma \ref{topnak}. We have an exact sequence of ${\mathcal
  L}(\widehat{N}^{(e)}_{0,1})$-modules$$0\longrightarrow D'({\mathcal
  F})\longrightarrow D({\mathcal F})\longrightarrow {\rm Hom}^{\rm
  ct}_{{\mathfrak o}}({\mathcal F}(\sigma),L/{\mathfrak
  o}).$$With ${\mathcal
  F}(x_+)^{{N}^{(e)}_0, p_L=0}$ also ${\mathcal F}(\sigma)$ and hence the cokernel of $D'({\mathcal
  F})\to D({\mathcal F})$ is finite.

Formula (\ref{classet1}) is immediately verified already by evaluating elements
of $D({\mathcal F})$ resp. $D'({\mathcal F})$ on the level of $0$-chains: it
simply follows from $tX_+^0\subset X_+^0$ (and $t^{-1}(tX_+^0)=X_+^0$). To see
formula (\ref{classet2}) we observe first that it follows from formula
(\ref{humphr}) that any element in $H_0(X_+,{\mathcal
  F})$ can be represented by a $0$-chain with support
in $$X_+^0-Z_+^{(m-1)}(x_+)=\bigcup_{n\in
  N_0^{(1)}/t^{m}N_0^{(1)}t^{-m}}nt^mX_+^0.$$But evaluating elements
of $D({\mathcal F})$ resp. $D'({\mathcal F})$ on such $0$-chains it is clear
that $$\sum_{n\in
  N_0^{(1)}/t^{m}N_0^{(1)}t^{-m}}n\varphi_{t^m}\psi_{t^m}n^{-1}$$ acts as the
identity on ${\mathcal O}^{(e)}_{{\mathcal E}}\otimes_{{\mathcal
  L}({N}^{(e)}_0)}D'({\mathcal F}){\cong} {\mathcal O}^{(e)}_{{\mathcal E}}\otimes_{{\mathcal
  L}({N}^{(e)}_0)}D({\mathcal F})$.\hfill$\Box$\\

We write ${\bf D}({\mathcal F})={\mathcal O}^{(e)}_{{\mathcal E}}\otimes_{{\mathcal
  L}({N}^{(e)}_0)}D'({\mathcal F}){\cong} {\mathcal O}^{(e)}_{{\mathcal E}}\otimes_{{\mathcal
  L}({N}^{(e)}_0)}D({\mathcal F})$. The key property exploited in Colmez' study of his functor ${\bf D}$ is that in the case $F={\mathbb
  Q}_p$ (and $e=1$) the module ${\bf D}({\mathcal F})$ is finitely generated over ${\mathcal O}^{(1)}_{{\mathcal E}}$. Namely:

\begin{satz}\label{colsa} Suppose $F={\mathbb
  Q}_p$ and $e=1$ and that the $k$-vector space ${\mathcal F}(x_+)^{{N}^{(1)}_0,
    p=0}$ has dimension $n<\infty$. Then ${\bf D}({\mathcal F})$ can be
  generated as an
  ${\mathcal O}^{(1)}_{{\mathcal E}}$-module by $n$ elements. 
\end{satz}

{\sc Proof:} The proof is parallel to that of Lemma \ref{topnak}, with Lemma
\ref{nauchinv} (in that proof) replaced by the following fact: for $F={\mathbb
  Q}_p$ and $e=1$ the natural map ${\mathcal F}(\sigma)\cong {\mathcal
  F}(x_+)^{{N}_0^{(1)}}\to H_0(X_+,{\mathcal F})^{{N}^{(1)}_0}$ is bijective;
this has e.g. been shown in \cite{jalg} Theorem 3.2 (notice that
$({N}^{(1)}_0,U_{x_+}^{(1)})=U_{\sigma}^{(1)}$). (In fact \cite{jalg} Theorem
3.2 is stated for the coefficient ring $k$ only, but the devissage arguments
used in the proof of \cite{jalg} Theorem 3.3 easily give the result in general.)\hfill$\Box$\\

 If $F={\mathbb Q}_p$ then ${\mathcal O}^{(1)}_{{\mathcal E}}$ is exactly the ring ${\mathcal O}_{{\mathcal E}}$ used in Fontaine's theory of $(\varphi,\Gamma)$-modules. Theorem \ref{colsa} describes Colmez' functor from finitely generated
${\mathfrak o}$-torsion representations $V$ of $G$ (generated by
$V^{U^{(1)}_{\sigma}}$), to finitely generated \'{e}tale
$(\varphi,\Gamma)$-modules. Namely, from $V$ pass to ${\mathcal F}_V^{(1)}\in {\mathfrak
  C}^{(1)}(\Lambda)$, then restrict the
$\widehat{P}^{(1)}_{\star}$-action on ${\bf D}({\mathcal F}_V^{(1)})$ to ${P}^{(1)}_{\star}$
and of the $(\widehat{P}^{(1)}_{\star})^{-1}$-action on ${\bf D}({\mathcal F}_V^{(1)})$
to $({P}^{(1)}_{\star})^{-1}$. 

Thus, we may regard Theorem \ref{fet} as some sort of variation, available for
arbitrary $F$ (and $e$), on Colmez' construction. For $F={\mathbb
  Q}_p$ and $e=1$ it is in fact an enhancement of Colmez' construction in that
it provides actions $\widehat{P}_{\star}^{(1)}$, reps. of $(\widehat{P}^{(1)}_{\star})^{-1}$, on ${D}'({\mathcal F}_V^{(1)})$ resp. on ${D}({\mathcal F}_V^{(1)})$ .

Let $F$ be a finite extension of ${\mathbb Q}_p$, and let us indicate the
dependence on $F$ in our definitions by the symbol $[F]$. Using the trace
map described in Lemma \ref{tracehatn} below one might try to pass from (finitely
generated, as in Theorem \ref{fet}) \'{e}tale
$((\widehat{P}_{\star}^{(e)})^{-1},\widehat{P}^{(e)}_{\star})$-pairs over
${\mathcal L}(\widehat{N}^{(e)}_{0,1})[F]$ to (finitely
generated) \'{e}tale
$((\widehat{P}_{\star}^{(e)})^{-1},\widehat{P}^{(e)}_{\star})$-pairs over
${\mathcal L}(\widehat{N}^{(e)}_{0,1})[{\mathbb Q}_p]$. From here, can one pass to suitable modules finitely generated even over ${\mathcal O}^{(e)}_{{\mathcal E}}[{\mathbb
  Q}_p]$ ?

\subsection{The pro-$p$-group $\widehat{N}_{0,1}^{(e)}$ and the Iwasawa algebra ${\mathcal L}(\widehat{N}_{0,1}^{(e)})$}

For $i\ge0$ consider the group homomorphism$$\epsilon_i:\prod_{a\in{\mathcal
    O}_F/{\mathfrak p}^i_F}{\mathfrak p}_F^{i+e-1}\longrightarrow{\rm
  Aut}_{\rm sets}({\mathcal O}_F)$$sending $(x_a)_{a\in {\mathcal
    O}_F/{\mathfrak p}^i_F}$ (with $x_a\in {\mathfrak p}_F^{i+e-1}$) to the
following permutation of (the set underlying) ${\mathcal O}_F$: it sends $x\in
{\mathcal O}_F$ to $x+x_{a(x)}$, where $a(x)\in {\mathcal O}_F/{\mathfrak
  p}^i_F$ is defined by $a(x)\equiv x$ modulo ${\mathfrak p}^i_F$. Next, for
$k\ge0$ define the map of sets$$\delta_k:(\prod_{a\in{\mathcal O}_F/{\mathfrak
    p}^k_F}{\mathfrak p}^{k+e-1}_F)\times\ldots\times(\prod_{a\in{\mathcal
    O}_F/{\mathfrak p}_F}{\mathfrak p}^{e}_F)\times{\mathfrak
  p}^{e-1}_F\longrightarrow {\rm Aut}_{\rm sets}({\mathcal
  O}_F),$$$$((x_{a,k})_{a\in{\mathcal O}_F/{\mathfrak
    p}^k_F},\ldots,(x_{a,1})_{a\in{\mathcal O}_F/{\mathfrak
    p}_F},x_{0,0})\mapsto\epsilon_k((x_{a,k})_{a\in{\mathcal
    O}_F/{\mathfrak
    p}^k_F})\circ\ldots\circ\epsilon_1((x_{a,1})_{a\in{\mathcal
    O}_F/{\mathfrak p}_F})\circ \epsilon_0(x_{0,0}).$$It induces a
map of sets$$\overline{\delta}_k:(\prod_{a\in{\mathcal O}_F/{\mathfrak
    p}^k_F}{\mathfrak p}^{k+e-1}_F)\times\ldots\times(\prod_{a\in{\mathcal
    O}_F/{\mathfrak p}_F}{\mathfrak p}^{e}_F)\times{\mathfrak
  p}^{e-1}_F\longrightarrow {\rm Aut}_{\rm sets}({\mathcal O}_F/{\mathfrak
  p}^{k+1}_F).$$Define $H^{(e)}_{0}={\mathfrak p}^{e-1}_F$. Let $H^{(e)}_{0}$
act on $\prod_{a\in{\mathcal O}_F/{\mathfrak p}_F}{\mathfrak p}_F^e$ by its
action $\overline{\delta}_0$ on the index set ${\mathcal O}_F/{\mathfrak
  p}_F$. Using this action we define $H^{(e)}_{1}=(\prod_{a\in{\mathcal
    O}_F/{\mathfrak p}_F}{\mathfrak p}_F^e)\rtimes H^{(e)}_{0}$. One easily
checks that with this group structure, $\overline{\delta}_1$ is a group
homomorphism. Using it we have $H^{(e)}_{1}$ acting on $\prod_{a\in{\mathcal
    O}_F/{\mathfrak p}^2_F}{\mathfrak p}^{e+1}_F$ and we may define
$H^{(e)}_{2}=(\prod_{a\in{\mathcal O}_F/{\mathfrak p}^2_F}{\mathfrak
  p}^{e+1}_F)\rtimes H^{(e)}_{1}$. Inductively we thus define a group
$H^{(e)}_{k}$ for every $k\ge0$. We then define $H^{(e)}_{\infty}$ to be the projective limit of the $H^{(e)}_{k}$, i.e.$$H^{(e)}_{\infty}=(\ldots((\prod_{a\in{\mathcal
    O}_F/{\mathfrak p}^k_F}{\mathfrak p}^{k+e-1}_F)\rtimes((\prod_{a\in{\mathcal
    O}_F/{\mathfrak p}^{k-1}_F}{\mathfrak
  p}^{k+e-2}_F)\rtimes(\ldots((\prod_{a\in{\mathcal O}_F/{\mathfrak
    p}_F}{\mathfrak p}_F^e)\rtimes{\mathfrak p}_F^{e-1})\ldots)))\ldots).$$

We have $H^{(e')}_{\infty}\subset H^{(e)}_{\infty}$ for $e'\ge e$ while, of
course, abstractly the $H^{(e)}_{\infty}$ for the various $e$ are all
isomorphic.

Define the group homomorphism $\partial^{(e)}:H^{(e+1)}_{\infty}\to H_{\infty}^{(e)}$ not as the inclusion, but as the product of all the maps $$({\rm diag},-{\rm id}):\prod_{a\in{\mathcal O}_F/{\mathfrak p}^k_F}{\mathfrak p}^{k+e}_F\longrightarrow(\prod_{a\in{\mathcal O}_F/{\mathfrak p}^{k+1}_F}{\mathfrak p}^{k+e}_F)\rtimes(\prod_{a\in{\mathcal O}_F/{\mathfrak p}^k_F}{\mathfrak p}^{k+e-1}_F),$$i.e. the negative of the inclusion $\prod_{a\in{\mathcal O}_F/{\mathfrak p}^k_F}{\mathfrak p}^{k+e}_F\subset \prod_{a\in{\mathcal O}_F/{\mathfrak p}^k_F}{\mathfrak p}^{k+e-1}_F$ in the second component, and the diagonal embedding $\prod_{a\in{\mathcal O}_F/{\mathfrak p}^k_F}{\mathfrak p}^{k+e}_F\hookrightarrow \prod_{a\in{\mathcal O}_F/{\mathfrak p}^{k+1}_F}{\mathfrak p}^{k+e}_F$ (induced by the natural projection of index sets ${\mathcal O}_F/{\mathfrak p}^k_F\to{\mathcal O}_F/{\mathfrak p}^{k+1}_F$) in the first component.

\begin{lem}\label{propexp} We have a natural isomorphism \begin{gather}\frac{H^{(e)}_{\infty}}{{\rm im}(\partial^{(e)})}\cong\widehat{N}_{0,1}^{(e)}\label{neun}.\end{gather}
\end{lem}

{\sc Proof:} $N$ acts simply transitively on the set of ends of $X$ different
from $t^{-\infty}x_+$; similarly, $N_0^{(1)}$ acts simply transitively on the set of (infinite)
ends of $X_+$. Therefore, identifying $N_0^{(1)}$ with the additive group
$({\mathcal O}_F,+)$ and sending $f\in {\mathcal O}_F\cong N_0^{(1)}$ to the end $f(t^{\infty}x_+)$ is a bijection between ${\mathcal O}_F$
and the set of ends of $X_+$. The group $H^{(e)}_{\infty}$ acts on (the set underlying) ${\mathcal O}_F$
through $\delta_{\infty}=\lim_{\leftarrow}\delta_k$. Together we get an action
of $H^{(e)}_{\infty}$ on the set of ends of $X_+$. It induces an action of $H^{(e)}_{\infty}$ on
$X^0_+$ as follows. Given $x\in X^0_+$, choose an end $\alpha$ of $X_+$ with $x\in [x_+,\alpha]$. For $h\in H^{(e)}_{\infty}$ we then define $hx\in X^0_+$ as the unique vertex belonging to $[x_+,h\alpha]$ with $d(x_+,x)=d(x_+,hx)$; it is independent of the choice of $\alpha$. It
follows from the definitions that this $H^{(e)}_{\infty}$-action on $X^0_+$ induces a surjective homomorphism
$H^{(e)}_{\infty}\to \widehat{N}_{0,1}^{(e)}$ whose kernel is ${\rm im}(\partial^{(e)})$.\hfill$\Box$\\

Let us redescribe $\widehat{N}^{(e)}_{0,1}$ and ${\mathcal
  L}(\widehat{N}^{(e)}_{0,1})$ in the case $F={\mathbb Q}_p$, using the
isomorphism (\ref{neun}). For $k\ge0$ and $i\in {\mathbb Z}_p/p^k{\mathbb
  Z}_p$ consider the corresponding copy of $p^{k+e-1}{\mathbb Z}_p$ as a
factor of $H_{\infty}^{(e)}$. Let us denote by $e^{(k)}_i$ its topological
generator $p^{k+e-1}$. Then $H_{\infty}^{(e)}$ is the pro-$p$-completion of
the group generated by all the $e^{(k)}_i$, subject only to the following
relations:$$e^{(k)}_i\cdot e^{(\ell)}_j=e^{(\ell)}_j\cdot
e^{(k)}_{i+p^{\ell}}\quad\quad\quad\mbox{ if }k\ge\ell\mbox{ and }i\equiv
j\mbox{ mod }(p^{\ell}),$$$$e^{(k)}_i\cdot e^{(\ell)}_j=e^{(\ell)}_j\cdot
e^{(k)}_{i}\quad\quad\quad\mbox{ if }k\ge\ell\mbox{ and }i\not\equiv j\mbox{
  mod }(p^{\ell}).$$(Here, of course, $i+p^{\ell}$ is the class of
$i+p^{\ell}$ modulo $(p^k)$.) Next, ${\rm im}(\partial^{(e)})$ is
topologically generated inside $H_{\infty}^{(e)}$  by all
expressions $$\prod_{i\in {\mathbb Z}_p/p^k{\mathbb Z}_p\atop i-j\in
  p^{k-1}{\mathbb Z}_p}e_i^{(k)}=(e_j^{(k-1)})^p$$for $k\ge1$ and
$j\in{\mathbb Z}_p/p^{k-1}{\mathbb Z}_p$. Write $U_i^{(k)}=1-e^{(k)}_i$ in
${\mathcal L}(H_{\infty}^{(e)})$ for $k\ge0$ and $i\in{\mathbb
  Z}_p/p^k{\mathbb Z}_p$. The above description of $H_{\infty}^{(e)}$
translates as follows: ${\mathcal L}(H_{\infty}^{(e)})$ is the ($p$-adically complete) formal power series algebra in the variables $U_i^{(k)}$, subject to the commutation relations$$U^{(k)}_i\cdot U^{(\ell)}_j=U^{(\ell)}_j\cdot U^{(k)}_{i+p^{\ell}}\quad\quad\quad\mbox{ if }k\ge\ell\mbox{ and }i\equiv j\mbox{ mod }(p^{\ell}),$$$$U^{(k)}_i\cdot U^{(\ell)}_j=U^{(\ell)}_j\cdot U^{(k)}_{i}\quad\quad\quad\mbox{ if }k\ge\ell\mbox{ and }i\not\equiv j\mbox{ mod }(p^{\ell}).$$Next, the quotient ${\mathcal L}(\widehat{N}^{(e)}_{0,1})$ is then obtained by dividing out the closure of the ideal generated by all expressions $$\prod_{i\in {\mathbb Z}_p/p^k{\mathbb Z}_p\atop i-j\in p^{k-1}{\mathbb Z}_p}(1-U^{(k)}_i)-(1-U_j^{(k-1)})^p$$for $k\ge1$ and $j\in{\mathbb Z}_p/p^{k-1}{\mathbb Z}_p$.\\

To indicate the dependence on the local field $F$ (which was fixed so far) let us now write $\widehat{N}^{(e)}_{0,1}[F]$ and $H_{\infty}^{(e)}[F]$ and $\partial^{(e)}[F]$ instead of $\widehat{N}^{(e)}_{0,1}$ and $H_{\infty}^{(e)}$ and $\partial^{(e)}$. Let $F/E$ be a finite extension of local fields, let $e(F/E)$ be its ramification index.

\begin{lem}\label{tracehatn} For any $e, e'\in{\mathbb N}$ with $e'\cdot e(F/E)\le e$ there is a natural trace map $${\rm tr}:\widehat{N}^{(e)}_{0,1}[F]\longrightarrow \widehat{N}^{(e')}_{0,1}[E]$$ 
\end{lem}

{\sc Proof:} Let ${\rm trace}_{F/E}:F\to E$ be
  the usual trace map. As $e'\cdot e(F/E)\le e$ it restricts to a map ${\rm trace}_{F/E}:{\mathfrak p}_F^{e+m\cdot e(F/E)-t-1}\to {\mathfrak p}_E^{e'+m-1}$ for all $m\ge0$ and all $0\le t\le e(F/E)-1$. For all $m\ge1$ and  $0\le t\le e(F/E)-1$ the inclusion ${\mathcal O}_E\subset {\mathcal O}_F$ induces an injective map ${\mathcal O}_E/{\mathfrak p}_E^{m}\to {\mathcal O}_F/{\mathfrak p}_F^{m\cdot e(F/E)-t}$. We may thus define the composition$$\prod_{{\mathcal O}_F/{\mathfrak p}_F^{m\cdot e(F/E)-t}}{\mathfrak p}_F^{e+m\cdot e(F/E)-t-1}\longrightarrow\prod_{{\mathcal O}_F/{\mathfrak p}_F^{m\cdot e(F/E)-t}}{\mathfrak p}_E^{e'+m-1}\longrightarrow \prod_{{\mathcal O}_E/{\mathfrak p}_E^{m}}{\mathfrak p}_E^{e'+m-1}$$where the first arrow is the product of the above trace maps (in each factor), and the second arrow is the natural projection induced by the above inclusion of index sets. Summing over all $0\le t\le e(F/E)-1$ we obtain a map\begin{gather}\prod_{{\mathcal O}_F/{\mathfrak p}_F^{m\cdot e(F/E)}}{\mathfrak p}_F^{e-1+m\cdot e(F/E)}\times\cdots\times \prod_{{\mathcal O}_F/{\mathfrak p}_F^{(m-1)\cdot e(F/E)+1}}{\mathfrak p}_F^{e+(m-1)\cdot e(F/E)}\longrightarrow\prod_{{\mathcal O}_E/{\mathfrak p}_E^{m}}{\mathfrak p}_E^{e'-1+m}\label{prodfacttra}\end{gather}for any $m\ge1$. Taking the product of the maps (\ref{prodfacttra}) for all $m\ge1$, while taking for $m=0$ just the trace map ${\rm trace}_{F/E}:{\mathfrak p}_F^{e-1}\to {\mathfrak p}_E^{e'-1}$, we obtain a trace map $H_{\infty}^{(e)}[F]\to H_{\infty}^{(e')}[E]$. It is straightforward to check that it maps ${\rm im}(\partial^{(e)}[F])$ into ${\rm im}(\partial^{(e)}[E])$, hence it induces, using the isomorphisms (\ref{neun}) for $\widehat{N}^{(e)}_{0,1}[F]$ and $\widehat{N}^{(e')}_{0,1}[E]$, a trace map ${\rm tr}:\widehat{N}^{(e)}_{0,1}[F]\longrightarrow \widehat{N}^{(e')}_{0,1}[E]$ as desired.\hfill$\Box$

\begin{flushleft} \textsc{Humboldt-Universit\"at zu Berlin\\Institut f\"ur Mathematik\\Rudower Chaussee 25\\12489 Berlin, Germany}\\ \textit{E-mail address}:
gkloenne@math.hu-berlin.de \end{flushleft} \end{document}